\documentclass[preprint,10.5pt]{elsarticle}

\usepackage{amssymb}
\usepackage{amsmath}

\usepackage{mathrsfs}
\usepackage{amsfonts}
\usepackage{CJK,fancybox,color,graphicx,amsmath,amsthm,amssymb,enumerate}

\newtheorem{Lemma}{Lemma}[section]

\newtheorem{thm}{Theorem}[section]

\theoremstyle{definition}

\theoremstyle{remark}

\numberwithin{equation}{section}

\newcommand{\D}{\displaystyle}
\newcommand{\DF}[2]{\frac{\D#1}{\D#2}}

\begin{document}

\begin{frontmatter}



\title{Concentration in the flux approximation limit of Riemann solutions to the extended Chaplygin gas equations with Coulomb-like friction \tnoteref{1}}

\tnotetext[1]{Supported by NSFC (71601085).}

\author{Qingling Zhang \corref{cor1}}
\cortext[cor1]{Corresponding author:zhangqingling2002@163.com}
\address{School of Mathematics and Computer Sciences, Jianghan University, Wuhan 430056, PR China }
\begin{abstract}

In this paper, two kinds of occurrence mechanism on the phenomenon
of concentration and the formation of delta shock waves are analyzed and identified in the flux
approximation limit of Riemann solutions to the extended Chaplygin
gas equations with Coulomb-like friction, whose special case can also be seen as the model of the magnetogasdynamics with Coulomb-like friction.
Firstly, by introducing a transformation, the Riemann
problem for the extended Chaplygin gas equations with Coulomb-like friction is solved
completely. Secondly, we rigorously show that, as the pressure
vanishes, any two-shock Riemann solution to the nonhomogeneous extended Chaplygin
gas equations tends to a $\delta$-shock solution to the correspongding nonhomogeneous transportation
equations, and the intermediate density between the two shocks tends
to a weighted $\delta$-measure that forms the $\delta$-shock; any
two-rarefaction-wave Riemann solution to the nonhomogeneous extended Chaplygin gas
equations tends to a two-contact-discontinuity solution to the correspongding nonhomogeneous
transportation equations, and the nonvacuum intermediate state between
the two rarefaction waves tends to a vacuum state. At last, we also
show that, as the pressure approaches the generalized Chaplygin
pressure, any two-shock Riemann solution to the nonhomogeneous extended Chaplygin
gas equations tends to a delta-shock
solution to the correspongding nonhomogeneous generalized Chaplygin gas equations. In a word, we have generalized all the
results about the vanishing pressure limit now available for homogeneous equations to the nonhomogeneous case.

\end{abstract}

\begin{keyword}
Extended Chaplygin gas; Delta shock waves; flux approximation limit;
Riemann solutions; transportation equations; generalized Chaplygin gas; Coulomb-like friction.

\MSC[2008] 35L65 \sep 35L67  \sep 35B25



\end{keyword}

\end{frontmatter}



\section{Introduction }

\setcounter{equation}{0}

The extended Chaplygin gas equations with Coulomb-like friction can be expressed as
\begin{equation}\label{1.1}
\left\{\begin{array}{ll}
\rho_t+(\rho u)_x=0,\\
(\rho u)_t+(\rho u^2+P)_x=\beta \rho,
\end{array}\right.
\end{equation}
where $\rho$, $u$ and $P$ represent the density, the velocity and
the scalar pressure, respectively, $\beta$ is a constant, and
\begin{equation}\label{1.2}
P=A\rho^{n}-\frac{B}{\rho^{\alpha}},\ \ \ 1 \leq n \leq 3,\ \
0<\alpha\leq1,
\end{equation}
with two parameters $A,B>0$. For $n=2$, this model can also be seen as the magnetogasdynamics with generalized Chaplygin pressure and Coulomb-like friction.

When $\beta=0$ in $(\ref{1.1})$, it becomes the extended Chaplygin gas equations
which was proposed by Naji in 2014
\cite{Naji} to study the evolution of dark energy. Moreover, it contains the magnetogasdynamics with generalized Chaplygin pressure as a special model with $n=2$.
When $B=0$ in
$(\ref{1.2})$, $P=A\rho^{n}$ is the standard equation of state of
perfect fluid. Up to now, various kinds of theoretical models have
been proposed to interpret the behavior of dark energy. Specially,
when $n=1$ in $(\ref{1.2})$, it reduces to the state equation for
modified Chaplygin gas, which was originally proposed by Benaoum in
2002 \cite{Benaoum}. As an exotic fluid, such a gas can explain the
current accelerated expansion of the universe. Whereas when $A=0$ in
$(\ref{1.2})$, $P=-\frac{B}{\rho^{\alpha}}$ is called the pressure
for the generalized Chaplygin gas \cite{Setare}. Furthermore, when
$\alpha=1$, $P=-\frac{B}{\rho}$ is called the pressure for (pure)
Chaplygin gas which was introduced by Chaplyin \cite{Chaplygin},
Tsien \cite{Tsien} and von Karman \cite{von Karman}. Such a gas owns a negative
pressure and occurs in certain theories of cosmology. It has also been
 advertised as a possible model for dark energy
\cite{Bilic-Tupper-Viollier,Gorini-Kamenshchik-Moschella-Pasquier}.

When two parameters $A$, $B\rightarrow0$, the limit system of
$(\ref{1.1})$ with $(\ref{1.2})$ formally becomes the following
transportation equations with Coulomb-like friction:
\begin{equation}\label{1.3}
\left\{\begin{array}{ll}
\rho_t+(\rho u)_x=0,\\
(\rho u)_t+(\rho u^2)_x=\beta \rho,
\end{array}\right.
\end{equation}
whose Riemann problem was firstly studied by Shen in \cite{Shen1}, which showed that the $\delta$-shock and
vacuum states do occur in Riemann solutions.

For $\beta=0$ in $(\ref{1.3})$, it becomes the transportation equations, which was also called the zero-pressure gas dynamics \cite{Bouchut,Brenier-Grenier,Li-Cao}. It can be used to describe some important physical phenomena, such as the motion of
free particles sticking together under collision and the formation
of large scale structures in the universe
\cite{E-Rykov-Sinai}. The transportation equations
have been studied extensively since 1994. The existence of measure solutions of the Riemann
problem was first proved by Bouchut \cite{Bouchut} and the existence
of the global weak solutions was obtained by Brenier and Grenier
\cite{Bouchut} and E.Rykov and Sinai \cite{E-Rykov-Sinai}. Sheng and
Zhang \cite{Sheng-Zhang} discovered that the $\delta$-shock and
vacuum states do occur in the Riemann solutions to the transport
equation by the vanishing viscosity method. For more results, one can refer to
\cite{Huang-Wang,Wang-Ding,Wang-Huang-Ding}.

$\delta$-shock is a kind of nonclassical nonlinear waves on which at
least one of the state variables becomes a singular measure.
Korchinski \cite{Korchinski} firstly introduced the concept of the
$\delta$-function into the classical weak solution in his
unpublished Ph. D. thesis. Tan, Zhang and Zheng
\cite{Tan-Zhang-Zheng} considered some 1-D reduced system and
discovered that the form of $\delta$-functions supported on shocks
was used as parts in their Riemann solutions for certain initial
data. LeFloch et al.\cite{LeFloch-Liu} applied the approach of
nonconservative product to consider nonlinear hyperbolic systems in
the nonconservative form. Recently, the weak asymptotic method was widely used to
study the $\delta$-shock wave type solution by Danilov and
Shelkovich et
al.\cite{Danilvo-Shelkovich1,Danilvo-Shelkovich2,Shelkovich}.

As for delta shock waves, much research is focused on exploring the
phenomena of concentration and and the formation of delta
shock waves in Riemann solutions. In \cite{Chen-Liu1},
Chen and Liu considered the Euler equations for isentropic fluids,
i.e., in $(\ref{1.1})$, they took the prototypical pressure function
as follows:
\begin{equation}\label{1.4}
P=\varepsilon\frac{\rho^\gamma}{\gamma},\ \ \gamma>1.
\end{equation}
They analyzed and identified the phenomena of concentration and
cavitation and the formation of $\delta$-shocks and vacuum states as
$\varepsilon\rightarrow 0$, which checked the numerical observation
for the 2-D case by Chang, Chen and Yang
\cite{Chang-Chen-Yang1,Chang-Chen-Yang2}. They also pointed out that
the occurrence of $\delta$-shocks and vacuum states in the process
of vanishing pressure limit can be regarded as a phenomenon of
resonance between the two characteristic fields. Moreover, they generalized this result
to the nonisentropic fluids in
\cite{Chen-Liu2}. Besides, the results were extended to the relativistic
Euler equations for polytropic gases in
\cite{Yin-Sheng}, the perturbed Aw-Rascle model in
\cite{Shen-Sun}, the magnetogasdynamics with generalized Chaplygin pressure in \cite{Chen-Sheng}, the modified Chaplygin gas equations \cite{Yang-Wang,Yang-Wang1}, and the nonhomogeneous (generalized) Chaplygin gas equations in \cite{Guo-Li-Yin1,Guo-Li-Yin2}, etc.

In this paper, in
contrast to the previous works in
\cite{Chen-Liu1,Chen-Liu2,Chen-Sheng,Shen-Sun,Yang-Wang,Yang-Wang1,Yin-Sheng} which concentrated on the homogeneous equations, we focus on the extended Chaplygin gas equations with Coulomb-like friction to discuss the phenomena of concentration and
cavitation and the formation of delta shock waves and vacuum states
in Riemann solutions as the pressure vanishes, or
tends to the generalized Chaplygin pressure.

It is noticed that, When $A$, $B\rightarrow0$, the system
$(\ref{1.1})$ with $(\ref{1.2})$ formally becomes the transportation
equations with Coulomb-like friction $(\ref{1.3})$. For fixed $B$, When $A\rightarrow0$, the
system $(\ref{1.1})$ with $(\ref{1.2})$ formally becomes the following
generalized Chaplygin gas equations with Coulomb-like friction
\begin{equation}\label{1.5}
\left\{\begin{array}{ll}
\rho_t+(\rho u)_x=0,\\
(\rho u)_t+(\rho u^2-\frac{B}{\rho^{\alpha}})_x=\beta\rho,
\end{array}\right.
\end{equation}
which corresponds to the Chaplygin gas equations with Coulomb-like friction when $\alpha=1$.
Recently, the research by Shen and Sun \cite{Shen2,Sun} showed that the $\delta$-shocks do occur in Riemann solutions to
(generalized) Chaplygin gas equations with Coulomb-like friction, but vacuum states do not occur. For more research about homogeneous (generalized) Chaplygin gas equations, one can refer to \cite{Brenier,Guo-Sheng-Zhang,Qu-Wang,Wang,Wang-Zhang2}.

In this paper, we first solve the Riemann problem of system
(\ref{1.1})-(\ref{1.2}) with Riemann initial data
\begin{equation}\label{1.6}
(\rho,u)(x,0)=(\rho_\pm,u_\pm),\ \ \ \ \pm x>0,
\end{equation}
where $\rho_\pm>0,\ u_\pm$ are arbitrary constants. By introducing a transformation, we reformulated the nonhomogeous equations into a conservative system.
With the help of the phase plane
analysis method, we constructed the Riemann solutions to (\ref{1.1})-(\ref{1.2}) and (\ref{1.6})
with four different structures: $R_1R_2$, $R_1S_2$, $S_1R_2$ and
$S_1S_2$.

Then we analyze the formation of $\delta$-shocks and vacuum states
in Riemann solutions as the pressure vanishes. It is shown that,
as the pressure vanishes, any two-shock Riemann solution to the
extended Chaplygin gas equations with Coulomb-like friction tends to a $\delta$-shock solution
to the transportation equations with Coulomb-like friction, and the intermediate density between the
two shocks tends to a weighted $\delta$-measure that forms the
$\delta$-shock; by contrast, any two-rarefaction-wave Riemann
solution to the extended Chaplygin gas equations with Coulomb-like friction tends to a
two-contact-discontinuity solution to the transportation equations with Coulomb-like friction, and
the nonvacuum intermediate state between the two rarefaction waves
tends to a vacuum state, even when the initial data stay away from
the vacuum. As a result, the delta shocks for the transportation
equations with Coulomb-like friction result from a phenomenon of concentration, while the
vacuum states results from a phenomenon of cavitation in the
vanishing pressure limit process. This is
consistent with those results obtained for homogeneous equations in
\cite{Chen-Liu1,Chen-Liu2,Chen-Sheng,Shen-Sun,Yang-Wang,Yang-Wang1,Yin-Sheng} and for nonhomogeneous equations in
\cite{Guo-Li-Yin1,Guo-Li-Yin2}.

In addition, we also proved that as the pressure tends to the
generalized Chaplygin pressure ($A\rightarrow0$), any two-shock
Riemann solution to the extended Chaplygin gas equations with Coulomb-like friction tends to a
$\delta$-shock solution to the generalized Chaplygin gas equations with Coulomb-like friction,
and the intermediate density between the two shocks tends to a
weighted $\delta$-measure that forms the $\delta$-shock.
Consequently, the delta shocks for the generalized Chaplygin gas
equations result from a phenomenon of concentration in the partly
vanishing pressure limit process. And we first generalized those results obtained for homogeneous equations in
\cite{Chen-Sheng,Yang-Wang1} to nonhomogeneous equations, which is also the novelty of this article lies in.

From the above analysis, we can find two kinds of occurrence
mechanism on the phenomenon of concentration and the formation of
delta shock wave in the flux approximation limit of the
extended Chaplygin gas equations with Coulomb-like friction. On one hand, since the strict hyperbolicity of the
limiting system $(\ref{1.3})$ fails, the delta shock
wave forms in the limit process as the pressure vanishes, see Section 5.
On the other hand, although the strict hyperbolicity of the limiting system
$(\ref{1.5})$ is preserved, the formation of delta
shock waves still occur as the pressure partly vanishes, see Section 6. In this
regard, the later occurrence
mechanism is much different from the former.
However, in any case, the phenomenon of concentration and the formation of
delta shock wave for the nonhomogeneous equations can be regarded as a process of resonance formation
between two characteristic fields, which well preserves and generalizes the results obtained for the homogeneous equations.

The paper is organized as follows. In Section 2, we restate the
Riemann solutions to transportation equations with Coulomb-like friction $(\ref{1.3})$ and the
generalized Chaplygin gas equations with Coulomb-like friction $(\ref{1.5})$. In Section 3, we
investigate the Riemann problem for a conservative system reformulated from the extended Chaplygin gas
equations with Coulomb-like friction $(\ref{1.1})$-$(\ref{1.2})$ and examine the dependence of
the Riemann solutions on the two parameters $A,B>0$. In Section 4, we obtain the Riemann solutions to the nonhomogeneous equations
$(\ref{1.1})$-$(\ref{1.2})$ with $(\ref{1.6})$. In Section 5,
we analyze the formation of delta shocks in the limit process of Riemann solutions to the extended Chaplygin
gas equations $(\ref{1.1})$-$(\ref{1.2})$ with $(\ref{1.6})$ as the
pressure vanishes. In Section 6, we discuss the formation of delta shocks in the limit process of Riemann
solutions to the extended Chaplygin gas equations
$(\ref{1.1})$-$(\ref{1.2})$ with $(\ref{1.6})$ as the pressure
approaches to the generalized Chaplygin pressure. Finally,
conclusions and discussions are drawn in Section 7.

\section{Preliminaries}
\subsection{Riemann solutions to the transportation equations with Coulomb-like friction}

In this
section, we restate the Riemann solutions to the nonhomogeneous transportation equations
$(\ref{1.3})$ with initial data $(\ref{1.6})$. See
\cite{Shen1} for more details.

The transportation equations $(\ref{1.3})$ have a double eigenvalue
$\lambda=u$ and only one right eigenvectors $ \vec{r}=(1,0)^T $.
Furthermore, we have $\nabla\lambda\cdot\vec{r}=0$, which means that
$\lambda$ is linearly degenerate. The Riemann problem $(\ref{1.3})$
and $(\ref{1.6})$ can be solved by contact discontinuities, vacuum
or $\delta$-shocks connecting two states
$(\rho_\pm,u_\pm+\beta t)$.

It is noticed that the characteristic equations of the system $(\ref{1.3})$ is
\begin{equation}\label{2.1}
\frac{dx}{dt}=u,\ \ \frac{du}{dt}=\beta.
\end{equation}

For the Riemann problem $(\ref{1.3})$ and $(\ref{1.6})$, the characteristic curve starting from the initial point $(x_{0},0)$
and the value of $(\rho,u)$ along this characteristic curve before intersection can be expressed, respectively, as
$$x(t)=u_{-}t+\frac{1}{2}\beta t^{2}+x_{0},\ \ \rho=\rho_{-},\ \ u=u_{-}+\beta t, x_{0}<0,$$
and
$$x(t)=u_{+}t+\frac{1}{2}\beta t^{2}+x_{0},\ \ \rho=\rho_{+},\ \ u=u_{+}+\beta t, x_{0}>0.$$

For $u_-<u_+$, there is no characteristic passing through
the region $\{x: u_{-}t+\frac{1}{2}\beta t^{2}<x<u_{+}t+\frac{1}{2}\beta t^{2}\}$, so the vacuum should appear in
the region. The solution can be expressed as
\begin{equation}\label{2.2}
(\rho,u)(x,t)=\left\{\begin{array}{ll}
(\rho_-,u_-+\beta t),\ \ \ \ -\infty<x<u_{-}t+\frac{1}{2}\beta t^{2},\\
vacuum,\ \ \ \ \ \ \ \ \ \ u_{-}t+\frac{1}{2}\beta t^{2}<x<u_{+}t+\frac{1}{2}\beta t^{2},\\
(\rho_+,u_++\beta t),\ \ \ \ u_{+}t+\frac{1}{2}\beta t^{2}<x<\infty.
\end{array}\right.
\end{equation}

For $u_-=u_+$, it is easy to see that the two states
$(\rho_\pm,u_\pm+\beta t)$ can be connected by a contact discontinuity $x(t)=u_{\pm}t+\frac{1}{2}\beta t^{2}$.

For $u_->u_+$, a solution containing a weighted
$\delta$-measure supported on a curve will be constructed. Let
$x=x(t)$ be a discontinuity curve, we consider a piecewise smooth
solution of $(\ref{1.3})$ in the form
\begin{equation}\label{2.3}
(\rho,u)(x,t)=\left\{\begin{array}{ll}
(\rho_-,u_-+\beta t),\ \ \ \ \ \ \ \ \ \ \ \ \ \ \ \ \ \ \ \ \ \ \ \ x<x(t),\\
(w(t)\delta(x-x(t)),u_\delta(t)),\ \ \ \ \ \ \ \ \ \ \ \ x=x(t),\\
(\rho_+,u_++\beta t),\ \ \ \ \ \ \ \ \ \ \ \ \ \ \ \ \ \ \ \ \ \ \ \ x>x(t),
\end{array}\right.
\end{equation}
in which $w(t)$ is the weight of the delta shock wave, $u_\delta(t)$ is the value of $u$
on this delta shock wave, and $u_\delta(t)-\beta t$ is constant on this delta shock wave.

To define the measure solutions, a two-dimensional weighted
$\delta$-measure $p(s)\delta_S$ supported on a smooth curve
$S=\{(x(s),t(s)):a<s<b\}$ can be defined as
\begin{equation}\label{2.4}
\langle
p(s)\delta_S,\psi(x(s),t(s))\rangle=\int_a^bp(s)\psi(x(s),t(s))\sqrt{{x'(s)}^2+{t'(s)}^2}ds,
\end{equation}
for any $\psi\in C_0^\infty(R\times R_{+})$.

For convenience, we usually select the parameter $s=t$ and use
$w(t)=\sqrt{1+{x'(t)}^2}p(t)$ to denote the strength of delta shock
wave from now on.

As shown in \cite{Shen1}, for any $\psi\in C_0^\infty(R\times
R_{+})$, the $\delta$-measure solution $(\ref{2.3})$ constructed
above satisfies
\begin{equation}\label{2.5}
\left\{\begin{array}{ll} \langle\rho,\psi_t\rangle+\langle\rho
u,\psi_x\rangle=0,
\\
\langle\rho u,\psi_t\rangle+\langle\rho u^2,\psi_x\rangle=-\langle\beta\rho,\psi\rangle,
\end{array}\right.
\end{equation}
 in which
$$\langle\rho,\psi\rangle=\int_0^\infty\int_{-\infty}^\infty\rho_0\psi dxdt+\langle w(t)\delta_S,\psi\rangle,$$
$$\langle\rho u,\psi\rangle=\int_0^\infty\int_{-\infty}^\infty\rho_0 u_0\psi dxdt+\langle w(t)u_\delta(t)\delta_S,\psi\rangle,$$
where
$$\rho_0=\rho_-+(\rho_+-\rho_-)H(x-x(t)),$$
$$\rho_0u_{0}=\rho_-(u_{-}+\beta t)+(\rho_+(u_{+}+\beta t)-\rho_-(u_{-}+\beta t))H(x-x(t)).$$

Substituting $(\ref{2.3})$ into $(\ref{2.5})$, one can derive the
generalized Rankine-Hugoniot conditions
\begin{equation}\label{2.6}
\left\{\begin{array}{ll}
\DF{dx(t)}{dt}=\sigma(t)=u_\delta(t),\\
\DF{dw(t)}{dt}=\sigma(t)[\rho]-[\rho u],\\
\DF{d(w(t) u_ \delta(t))}{dt}=\sigma(t)[\rho u]-[\rho u^2]+\beta w(t),
\end{array}\right.
\end{equation}
where $[\rho]=\rho_+-\rho_-, [\rho u]=\rho_+(u_{+}+\beta t)-\rho_-(u_{-}+\beta t), [\rho u^2]=\rho_+(u_{+}+\beta t)^2-\rho_-(u_{-}+\beta t)^2$, denotes the jump of the function $\rho$, $\rho u$ and $\rho u^2$ across the $\delta$-shock.

Through solving $(\ref{2.6})$ with $x(0)=0,\ w(t)=0$, under the following $\delta$-entropy condition:
\begin{equation}\nonumber
u_++\beta t<\sigma(t)<u_-+\beta t,
\end{equation}we obtain
\begin{equation}\label{2.7}
\left\{\begin{array}{ll}
u_\delta(t)=\sigma(t)=\sigma_0+\beta t,\\
x(t)=\sigma_{0} t+\frac{1}{2}\beta t^2,\\
w(t)=\sqrt{\rho_-\rho_+}(u_--u_+)t,
\end{array}\right.
\end{equation}
where $\sigma_0=\DF{\sqrt{\rho_-}u_-+\sqrt{\rho_+}u_+}{\sqrt{\rho_-}+\sqrt{\rho_+}}.$

\subsection{Delta shock waves for Riemann solutions to the generalized Chaplygin gas
equations with Coulomb-like friction}

In this section, for simplicity, we only restate the delta shock wave solution to the Riemann problem for the nonhomogenous generalized
Chaplygin gas equations $(\ref{1.5})$-$(\ref{1.6})$. For details about the Riemann problem, one
can refer to \cite{Shen2, Sun}.

It is easy to see that $(\ref{1.5})$ has two eigenvalues
$$\lambda_1^B=u-\sqrt{\alpha B}\rho^{-\frac{\alpha+1}{2}},\ \ \lambda_2^B=u+\sqrt{\alpha B}\rho^{-\frac{\alpha+1}{2}},$$
with corresponding right eigenvectors
$$\overrightarrow{r_1}^B=(-\sqrt{\alpha B}\rho^{-\frac{\alpha+1}{2}},\rho)^T,\ \ \overrightarrow{r_2}^B=(\sqrt{\alpha B}\rho^{-\frac{\alpha+1}{2}},\rho)^T.$$
So $(\ref{1.5})$ is strictly hyperbolic for $\rho>0$. Moreover, when
$0<\alpha<1$, we have
 $\bigtriangledown\lambda_i^B\cdot\overrightarrow{r_i}^B\neq0$,
$i=1,2$, which implies that $\lambda_1^B$ and $\lambda_2^B$ are both
genuinely nonlinear and the associated waves are rarefaction waves
and shock waves. When $\alpha=1$, $\bigtriangledown\lambda_i^B\cdot
\overrightarrow{r_i}^B=0$, $i=1,2$, which implies that $\lambda_1^B$
and $\lambda_2^B$ are both linearly degenerate and the associated
waves are both contact discontinuities, see \cite{Smoller}.

For given inital data $(\ref{1.6})$, when $u_++\sqrt{
B}\rho_+^{-\frac{\alpha+1}{2}} \leq u_--\sqrt{
B}\rho_-^{-\frac{\alpha+1}{2}}$, a delta shock wave must
develop in solutions. Under the definition (\ref{2.4}), a delta
shock wave can be introduced to construct the solution of
(\ref{1.5})-(\ref{1.6}), which can be expressed as
\begin{equation}\label{2.8}
(\rho,u)(x,t)=\left\{\begin{array}{ll}
(\rho_-,u_-+\beta t),\ \ \ \ \ \ \ \ \ \ \ \ \ \ \ \ \ \ \ \ \ \ \ x<x^B (t),\\
(w^B(t)\delta(x-x^B (t)),u_\delta^{B}(t)),\ \ \  \ \ x=x^B (t),\\
(\rho_+,u_++\beta t),\ \ \ \ \ \ \ \ \ \ \ \ \ \ \ \ \ \ \ \ \ \ \
x>x^B(t),
\end{array}\right.
\end{equation}
where $x^{B}(t)$, $w^{B}(t)$ and $u_\delta^{B}(t)$ are respectively
denote the location and weight of the
$\delta$-shock, $u_\delta(t)$ is the value of $u$
on this delta shock wave, and $u_\delta(t)-\beta t$ is constant on this delta shock wave.

 As show in \cite{Shen2, Sun}, for any $\psi\in C_0^\infty(R\times
R_{+})$, the $\delta$-measure solution $(\ref{2.8})$ constructed
above satisfies
\begin{equation}\label{2.9}
\left\{\begin{array}{ll} \langle\rho,\psi_t\rangle+\langle\rho
u,\psi_x\rangle=0,
\\
\langle\rho u,\psi_t\rangle+\langle\rho u^2-\frac{B}{\rho^{\alpha}},\psi_x\rangle=-\langle\beta\rho,\psi\rangle,
\end{array}\right.
\end{equation}
 in which
$$\langle\rho u^2-\frac{B}{\rho^{\alpha}},\psi\rangle=\int_0^\infty\int_{-\infty}^\infty(\rho_0 u_0^2-\frac{B}{\rho_0^{\alpha}})\psi dxdt+\langle w(t)(u_\delta(t))^2\delta_S,\psi\rangle,$$
where
$$\rho_0u_{0}^2=\rho_-(u_{-}+\beta t)^2+(\rho_+(u_{+}+\beta t)^2-\rho_-(u_{-}+\beta t)^2)H(x-x(t)).$$
Moreover, $\frac{B}{\rho^{\alpha}}$ is equal to zero on the delta shock wave.

By the above weak solution definition in Subsection 2.1, for the system
(\ref{1.5}) with the delta shock wave solution (\ref{2.8}), we can get the following generalized Rankine-Hugoniot
conditions
\begin{equation}\label{2.10}
\left\{\begin{array}{ll}
\DF{dx^{B}(t)}{dt}=u_\delta^{B}(t)=\sigma^B(t),\\[4pt]
\DF{dw^{B}(t)}{dt}=\sigma^B(t)[\rho]-[\rho u],\\[4pt]
\DF{d(w^{B}(t)u_\delta^{B}(t))}{dt}=\sigma^B(t)[\rho u]-[\rho
u^2-\DF{B}{\rho^{\alpha}}]+\beta w(t),
\end{array}\right.
\end{equation}
where $[\rho u^2-\frac{B}{\rho^{\alpha}}]=\rho_+(u_{+}+\beta t)^2-\frac{B}{\rho_+^{\alpha}}-(\rho_-(u_{-}+\beta t)^2-\frac{B}{\rho_-^{\alpha}})$.

Then by solving (\ref{2.10}) with initial data $x(0)=0,\ w^B(0)=0$,
under the entropy condition
\begin{equation}\label{2.11}
u_++\sqrt{\alpha B}\rho_+^{-\frac{\alpha+1}{2}}+\beta t<\sigma^B(t)<
u_--\sqrt{\alpha  B}\rho_-^{-\frac{\alpha+1}{2}}+\beta t,
\end{equation}
 we can obtain
\begin{equation}\label{2.12}
\sigma^B(t)=u_\delta^{B}(t)=\sigma_{0}^{B}+\beta t,
\end{equation}
\begin{equation}\label{2.13}
x^B(t)=\sigma_{0}^{B}t+\frac{1}{2}\beta t^2,
\end{equation}
\begin{equation}\label{2.14}
w^B(t)=w_0^Bt,
\end{equation}
when $\rho_+\neq\rho_-$, where
 \begin{equation}\label{2.15}
 \sigma_{0}^{B}=\frac{\rho_+
u_+-\rho_-u_-+w_0^B}{\rho_+-\rho_-},
\end{equation}
\begin{equation}\label{2.16}
w_0^B=\big\{\rho_+\rho_-\big((u_+-u_-)^2-(\frac{1}{\rho_+}-\frac{1}{\rho_-})(\frac{B}{\rho_+^\alpha}-\frac{B}{\rho_-^\alpha})\big)\big\}^\frac{1}{2},
\end{equation}
and
\begin{equation}\label{2.17}
\sigma^B(t)=u_\delta^{B}(t)=\frac{1}{2}(u_++u_-)+\beta t,
\end{equation}
\begin{equation}\label{2.18}
x^B(t)=\frac{1}{2}(u_++u_-)t++\frac{1}{2}\beta t^2,
\end{equation}
\begin{equation}\label{2.19}
w^B(t)=(\rho_-u_--\rho_+ u_+)t,
\end{equation}
when $\rho_+=\rho_-$.

\section{Riemann problem for a modified conservative system from (\ref{1.1})}

In this section, we are devoted to the study of the Riemann problem
for a conservative system from (\ref{1.1}) in detail. Let us introduce
the new velocity $v(x,t)=u(x,t)-\beta t$, then the system
(\ref{1.1}) can be reformulated into a conservative form as follows:
\begin{equation}\label{3.1}
\left\{\begin{array}{ll}
\rho_t+(\rho (v+\beta t))_x=0,\\
(\rho v)_t+(\rho (v(v+\beta t) +P)_x=0.
 \end{array}\right.
\end{equation}
In fact, the change of variable was introduced by Faccanoni and
Mangeney \cite{Faccanoni-Mangeney} to study the shock and
rarefaction waves of the Riemann problem for the shallow water
equations with a Coulomb-like friction. Here, we also use this
transformation to study the shock and
rarefaction waves of the Riemann problem for the system
(\ref{1.1}).

Now we want to deal with the Riemann problem for the conservative
system (\ref{3.1}) and (\ref{1.2}) with the same
 Riemann initial data (\ref{1.6}) as follows:
\begin{equation}\label{3.2}
(\rho,v)(x,0)=\left\{\begin{array}{ll}
(\rho_-,u_-),\ \ x<0,\\
(\rho_ +,u_+),\ \ x>0.
\end{array}\right.
\end{equation}
We shall see hereafter that the Riemann solutions to (\ref{1.1})-(\ref{1.2}) and
(\ref{1.6}) can be obtained immediately from the Riemann solutions
to (\ref{3.1}) and (\ref{3.2}) by using the transformation of state
variables $(\rho,u)(x,t)=(\rho,v+\beta t)(x,t)$.

The eigenvalues of the system $(\ref{3.1})$ and $(\ref{1.2})$ are
$$\lambda_1^{AB}(\rho,v)=v+\beta t-\sqrt{An\rho^{n-1}+\frac{\alpha B}{\rho^{\alpha+1}}}, \quad \lambda_2^{AB}(\rho,v)=v+\beta t+\sqrt{An\rho^{n-1}+\frac{\alpha B}{\rho^{\alpha+1}}},$$
with corresponding right eigenvectors
$$\vec{r}_1^{AB}=(-\rho,\sqrt{An\rho^{n-1}+\frac{\alpha B}{\rho^{\alpha+1}}})^T,\quad \vec{r}_2^{AB}=(\rho,\sqrt{An\rho^{n-1}+\frac{\alpha B}{\rho^{\alpha+1}}})^T.$$

Moreover, we have
$$\nabla\lambda_i^{AB}\cdot \vec{r}_i^{AB}=\DF{An(n+1)\rho^{n+\alpha}+(1-\alpha)\alpha B}{2\sqrt{(An\rho^{n+\alpha}+\alpha B)\rho^{\alpha+1}}}>0\ \ (i=1,2).$$
Thus $\lambda_1^{AB}$ and $\lambda_2^{AB}$ are genuinely nonlinear
and the associated elementary waves are
shock waves denoted by $S^{AB}$ or rarefaction waves denoted by $R^{AB}$.
For details about above elementaty waves, one can refer to \cite{Smoller} to see how to solve the Riemann problem.

According to
\cite{Smoller}, for a given left state $(\rho_-,v_-)$, the
rarefaction wave curves in the phase plane, which are the sets of
states that can be connected on the right by a 1-rarefaction wave or
2-rarefaction wave, are as follows:
\begin{equation}\label{3.3}
R_1^{AB}(\rho_-,v_-): \left\{\begin{array}{ll} \xi=\lambda_1=v+\beta t-\sqrt{An\rho^{n-1}+\frac{\alpha B}{\rho^{\alpha+1}}},\\
v-v_-=-\int_{\rho_{-}}^{\rho}\frac{\sqrt{An\rho^{n-1}+\frac{\alpha
B}{\rho^{\alpha+1}}}}{\rho}d\rho,
\end{array}\right.
\end{equation}
and
\begin{equation}\label{3.4}
R_2^{AB}(\rho_-,v_-): \left\{\begin{array}{ll} \xi=\lambda_2=v+\beta t+\sqrt{An\rho^{n-1}+\frac{\alpha B}{\rho^{\alpha+1}}},\\
v-v_-=\int_{\rho_{-}}^{\rho}\frac{\sqrt{An\rho^{n-1}+\frac{\alpha
B}{\rho^{\alpha+1}}}}{\rho}d\rho.
\end{array}\right.
\end{equation}

From $(\ref{3.3})$ and $(\ref{3.4})$, we obtain that

\begin{eqnarray}
\frac{d\lambda_1^{AB}}{d\rho}=\frac{\partial\lambda_1^{AB}}{\partial
v}\frac{dv}{d\rho}+\frac{\partial\lambda_1^{AB}}{\partial\rho}
=-\DF{An(n+1)\rho^{n-1}+\frac{\alpha(1-\alpha)
B}{\rho^{\alpha+1}}}{2\rho\sqrt{An\rho^{n-1}+\frac{\alpha
B}{\rho^{\alpha+1}}}}<0\label{3.5}
\end{eqnarray}

\begin{eqnarray}
\frac{d\lambda_2^{AB}}{d\rho}=\frac{\partial\lambda_2^{AB}}{\partial
v}\frac{dv}{d\rho}+\frac{\partial\lambda_2^{AB}}{\partial\rho}
=\DF{An(n+1)\rho^{n-1}+\frac{\alpha(1-\alpha)
B}{\rho^{\alpha+1}}}{2\rho\sqrt{An\rho^{n-1}+\frac{\alpha
B}{\rho^{\alpha+1}}}}>0\label{3.6}
\end{eqnarray}
which imply that the velocity of 1-rarefaction (2-rarefaction) wave
$\lambda_1^{AB}$ ($\lambda_2^{AB}$) is monotonic decreasing
(increasing) with respect to $\rho$.

With the requirement
$\lambda_1^{AB}(\rho_{-},v_{-})<\lambda_1^{AB}(\rho,v)$ and
$\lambda_2^{AB}(\rho_{-},v_{-})<\lambda_2^{AB}(\rho,v)$, noticing
$(\ref{3.5})$ and $(\ref{3.6})$, we get that

\begin{equation}\label{3.7}
R_1^{AB}(\rho_-,v_-): \left\{\begin{array}{ll} \xi=\lambda_1=v+\beta t-\sqrt{An\rho^{n-1}+\frac{\alpha B}{\rho^{\alpha+1}}},\\
v-v_-=-\int_{\rho_{-}}^{\rho}\frac{\sqrt{An\rho^{n-1}+\frac{\alpha
B}{\rho^{\alpha+1}}}}{\rho}d\rho,\ \  \rho<\rho_{-},
\end{array}\right.
\end{equation}
and
\begin{equation}\label{3.8}
R_2^{AB}(\rho_-,v_-): \left\{\begin{array}{ll} \xi=\lambda_2=v+\beta t+\sqrt{An\rho^{n-1}+\frac{\alpha B}{\rho^{\alpha+1}}},\\
v-v_-=\int_{\rho_{-}}^{\rho}\frac{\sqrt{An\rho^{n-1}+\frac{\alpha
B}{\rho^{\alpha+1}}}}{\rho}d\rho,\ \  \rho>\rho_{-}.
\end{array}\right.
\end{equation}

For the 1-rarefaction wave, through differentiating $v$ respect to
$\rho$ in the second equation in $(\ref{3.7})$, we get
\begin{equation}\label{3.9}
v_{\rho}=-\DF{\sqrt{An\rho^{n-1}+\frac{\alpha
B}{\rho^{\alpha+1}}}}{\rho}<0.
\end{equation}

\begin{equation}\label{3.10}
v_{\rho\rho}=\DF{-An(n-3)\rho^{n+\alpha}+\alpha(\alpha+3)
B}{2\rho^{2}\sqrt{An\rho^{n+\alpha}+\alpha B\rho^{\alpha+1}}}.
\end{equation}

Thus, it is easy to get $v_{\rho\rho}>0$ for $1\leq n\leq3$, i.e.,
the 1-rarefaction wave is convex for $1\leq n\leq3$ in the upper
half phase plane ($\rho>0$).

In addition, from the second equation of $(\ref{3.7})$, we have
$$v-v_-=\int_{\rho}^{\rho_{-}}\frac{\sqrt{An\rho^{n-1}+\frac{\alpha
B}{\rho^{\alpha+1}}}}{\rho}d\rho\geq\int_{\rho}^{\rho_{-}}\sqrt{\alpha
B}\rho^{-\frac{\alpha+1}{2}-1}d\rho=\frac{2\sqrt{\alpha
B}}{\alpha+1}(\rho^{-\frac{\alpha+1}{2}}-\rho_{-}^{-\frac{\alpha+1}{2}}),$$
which means that $\lim\limits_{\rho\rightarrow0}v=+\infty.$

By a similar computation, we have that, for the 2-rarefaction wave,
$v_{\rho}>0$, $v_{\rho\rho}<0$ for $1\leq n\leq3$ and
$\lim\limits_{\rho\rightarrow+\infty}v=+\infty.$ Thus, we can draw
the conclusion that the 2-rarefaction wave is concave for $1\leq
n\leq3$ in the upper half phase plane ($\rho>0$).

Let us return our attention on the shock wave which is a piecewise
constant discontinuous solution, satisfying the Rankine-Hugoniot
conditions and the entropy condition. Here  the Ranking-Hugoniot
conditions can be derived in a standard method as in \cite{Smoller},
since the parameter $t$ only appears
 in the flux functions in the conservative
system (\ref{3.1}). For a bounded discontinuity at $x=x(t)$, let us
denote $\sigma^{AB}(t)=x'(t)$, then the Ranking-Hugoniot conditions for
the conservative system (\ref{3.1}) can be expressed as
\begin{equation}\label{3.11}
\left\{\begin{array}{ll}
\sigma^{AB}(t)[\rho]=[\rho (v+\beta t)],\\
\sigma^{AB}(t)[\rho v]=[\rho v(v+\beta t)+P],\ \
P=A\rho^{n}-\frac{B}{\rho^{\alpha}},
\end{array}\right.
\end{equation}
where $[\rho]=\rho_{r}-\rho_l$ with
$\rho_{l}=\rho(x(t)-0,t)$,\ $\rho_{r}=\rho(x(t)+0,t)$, in which
$[\rho]$ denote the jump of $\rho$ across the discontinuity, etc. It
is clear that the propagation speed of the discontinuity depends on
the parameter $t$, which is obviously different from classical
hyperbolic conservation laws.

Eliminating $\sigma$ from $(\ref{3.11})$, we obtain
\begin{equation}\label{3.12}
v-v_{-}=\pm\sqrt{\frac{\rho-\rho_{-}}{\rho\rho_{-}}\Big(A(\rho^{n}-\rho_{-}^{n})-B(\frac{1}{\rho^{\alpha}}-\frac{1}{\rho_{-}^{\alpha}})\Big)}.
\end{equation}

Using the Lax entropy condition, the 1-shock satisfies
\begin{equation}\label{3.13}
\sigma^{AB}(t)<\lambda_{1}^{AB}(\rho_{-},v_{-}),\ \
\lambda_{1}^{AB}(\rho,v)<\sigma^{AB}(t)<\lambda_{2}^{AB}(\rho,v),
\end{equation}
while the 1-shock satisfies
\begin{equation}\label{3.14}
\lambda_{1}^{AB}(\rho_{-},v_{-})<\sigma^{AB}(t)<\lambda_{2}^{AB}(\rho_{-},v_{-}),\
\ \lambda_{2}^{AB}(\rho,v)<\sigma^{AB}(t).
\end{equation}

From the first equation in $(\ref{3.11})$, we have
\begin{equation}\label{3.15}
\sigma^{AB}(t)=\frac{\rho
(v+\beta t)-\rho_{-}v_{-}}{\rho-\rho_{-}}=v+\frac{\rho_{-}(
v-v_{-})}{\rho-\rho_{-}}+\beta t=v_{-} +\frac{\rho(
v-v_{-})}{\rho-\rho_{-}}+\beta t.
\end{equation}

Thus, by a simple calculation, $(\ref{3.13})$ is equivalent to
\begin{equation}\label{3.16}
-\rho\sqrt{An\rho^{n-1}+\frac{\alpha
B}{\rho^{\alpha+1}}}<\frac{\rho\rho_{-}(
v-v_{-})}{\rho-\rho_{-}}<-\rho_{-}\sqrt{An\rho_{-}^{n-1}+\frac{\alpha
B}{\rho_{-}^{\alpha+1}}},
\end{equation}
and $(\ref{3.14})$ is equivalent to
\begin{equation}\label{3.17}
\rho\sqrt{An\rho^{n-1}+\frac{\alpha
B}{\rho^{\alpha+1}}}<\frac{\rho\rho_{-}(
v-v_{-})}{\rho-\rho_{-}}<\rho_{-}\sqrt{An\rho_{-}^{n-1}+\frac{\alpha
B}{\rho_{-}^{\alpha+1}}}.
\end{equation}
$(\ref{3.14})$ and $(\ref{3.15})$ imply that $\rho>\rho_{-}$,
$u<u_{-}$ and $\rho<\rho_{-}$, $v<v_{-}$, respectively.

 Through the above analysis, for a given left
state $(\rho_-,v_-)$, the shock curves in the phase plane, which are
the sets of states that can be connected on the right by a 1-shock
or 2-shock, are as follows:

\
\begin{equation}\label{3.18}
S_1^{AB}(\rho_-,v_-): \left\{\begin{array}{ll}
\sigma_1=\DF{\rho v-\rho_-v_-}{\rho-\rho_-}+\beta t,\\
v-v_-=-\sqrt{\frac{\rho-\rho_{-}}{\rho\rho_{-}}\Big(A(\rho^{n}-\rho_{-}^{n})-B(\frac{1}{\rho^{\alpha}}-\frac{1}{\rho_{-}^{\alpha}})\Big)},\
\ \rho>\rho_-,
\end{array}\right.
\end{equation}
and
\begin{equation}\label{3.19}
S_2^{AB}(\rho_-,v_-): \left\{\begin{array}{ll}
\sigma_2=\DF{\rho v-\rho_-v_-}{\rho-\rho_-}+\beta t,\\
v-v_-=-\sqrt{\frac{\rho-\rho_{-}}{\rho\rho_{-}}\Big(A(\rho^{n}-\rho_{-}^{n})-B(\frac{1}{\rho^{\alpha}}-\frac{1}{\rho_{-}^{\alpha}})\Big)},\
\ \rho<\rho_-.
\end{array}\right.
\end{equation}

For the 1-shock wave, through differentiating $u$ respect to $\rho$
in the second equation in $(\ref{3.16})$, we get
\begin{equation}\label{3.20}
2(v-v_{-})v_{\rho}=\frac{1}{\rho^{2}}\Big(A(\rho^{n}-\rho_{-}^{n})-B(\frac{1}{\rho^{\alpha}}-\frac{1}{\rho_{-}^{\alpha}})\Big)+
\frac{\rho-\rho_{-}}{\rho\rho_{-}}(An\rho^{n-1}+\frac{\alpha
B}{\rho^{\alpha+1}})>0,
\end{equation}
which means that $v_{\rho}<0$ for the 1-shock wave and that the
1-shock wave curve is starlike with respect to $(\rho_-,v_-)$ in the
region $\rho>\rho_-$. Similarly, we can get  $v_{\rho}>0$ for the
2-shock wave and that the 2-shock wave curve is starlike with
respect to $(\rho_-,v_-)$ in the region $\rho<\rho_-$. In addition,
it is easy to check that
$\lim\limits_{\rho\rightarrow+\infty}v=-\infty$ for the 1-shock wave
and
 $\lim\limits_{\rho\rightarrow0}v=-\infty$ for the 2-shock wave.

Through the analysis above, for a given left state $(\rho_-,v_-)=(\rho_-,u_-)$,
the sets of states connected with $(\rho_-,v_-)$ on the right in the
phase plane consist of the 1-rarefaction wave curve
$R_1^{AB}(\rho_-,v_-)$, the 2-rarefaction wave curve $R_2^{AB}(\rho_-,v_-)$,
the 1-shock curve $S_1^{AB}(\rho_-,v_-)$ and the 2-shock curve
$S_2^{AB}(\rho_-,v_-)$. These curves divide the upper half plane into
four parts $\rm{I}$, $\rm{I\!I}$,
$\rm{I\!I\!I}$ and $\rm{I\!V}$. Now, we put
all of these curves together in the upper half plane ($\rho>0$,
$v\in R$) to obtain a picture as in Fig.1.

By the phase plane analysis method, it is easy to construct
Riemann solutions for any given right state $(\rho_+,v_+)=(\rho_+,u_+)$ as follows:\\
\ \ (1)  $(\rho_+,v_+)\in\rm{I}:R_{1}^{AB}+R_{2}^{AB};$ \ \
(2)  $(\rho_+,v_+)\in\rm{I\!I}:S_{1}^{AB}+R_{2}^{AB};$\\
\ \ (3)  $(\rho_+,v_+)\in\rm{I\!I\!I}:R_{1}^{AB}+S_{2}^{AB};$ \
\ \
 (4) $(\rho_+,v_+)\in\rm{I\!V}:S_{1}^{AB}+S_{2}^{AB}.$


\unitlength 1.00mm 
\linethickness{0.4pt}
\ifx\plotpoint\undefined\newsavebox{\plotpoint}\fi 
\begin{picture}(98.71,55.84)(15,0)
\put(131.71,11.89){\vector(1,0){.07}}
\put(52.22,11.89){\line(1,0){79.49}}
\qbezier(57.57,16.35)(92.95,19.92)(111.08,49.84)
\qbezier(120.58,14.94)(78.31,20.94)(71.54,50.84)
\put(49.1,47.87){\vector(0,1){.07}}
\put(49.1,24.68){\line(0,1){23.19}}
\put(133.34,11.89){\makebox(0,0)[cc]{$v$}}
\put(49.1,49.59){\makebox(0,0)[cc]{$\rho$}}%
 \put(88,26.30){\circle*{.9}}%
 \put(96.71,26.02){\makebox(0,0)[cc]{$\scriptstyle(\rho_-,u_-)$}}
 \put(110.49,30.47){\makebox(0,0)[cc]{$\rm{I}$}}
\put(87.3,42.96){\makebox(0,0)[cc]{$\rm{I\!I}$}}
\put(67.38,30.77){\makebox(0,0)[cc]{$\rm{I\!V}$}}
\put(84.92,16.2){\makebox(0,0)[cc]{$\rm{I\!I\!I}$}}
\put(68.76,47.12){\makebox(0,0)[cc]{$S_1^{AB}$}}
\put(112.99,46.08){\makebox(0,0)[cc]{$R_2^{AB}$}}
\put(61.14,19.58){\makebox(0,0)[cc]{$S_2^{AB}$}}
\put(110,19.43){\makebox(0,0)[cc]{$R_1^{AB}$}}
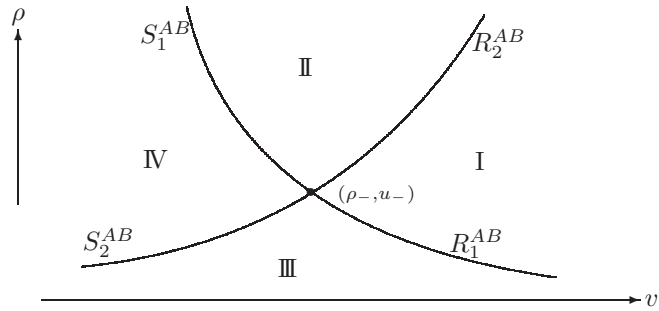
 \put(94.41,3.99){\makebox(0,0)[cc]{ Fig.1 the $(\rho,v)$ phase plane for the conservative system (\ref{3.1})}}
\end{picture}

\section{Riemann problem for the extended Chaplygin gas equations with Coulomb-like friction}

 In this section, let us return to the Riemann problem (\ref{1.1})-(\ref{1.2})
and (\ref{1.6}). If $(\rho_+,u_+)\in \rm{I}$ , the Riemann solutions to
(\ref{1.1})-(\ref{1.2}) and (\ref{1.6})  $R_1^{AB}+R_{2}^{AB}$ can be represented as
\begin{equation}\label{4.1}
(\rho,u)(x,t)=\left\{\begin{array}{ll}
(\rho_-,u_-+\beta t),\ \ \ \ \ \ \ \ \ \ \ x<x_{1}^-(t),\\
R_{1}^{AB},\ \ \  \ \ \ \ \ \ \ \ \ \ \ \ \ \ \ \ \ \ \ x_{1}^-(t)<x<x_{1}^+(t),\\
(\rho_*,v_*+\beta t),\ \ \  \ \ \ \ \ \ \ \ \ x_{1}^+(t)<x<x_{2}^-(t),\\
R_{2}^{AB},\ \ \  \ \ \ \ \ \ \ \ \ \ \ \ \  \ \ \ \ \ \ x_{2}^-(t)<x<x_{2}^+(t),\\
(\rho_+,u_++\beta t),\ \ \ \ \ \ \ \ \ \ \ x>x_{2}^+(t),
\end{array}\right.
\end{equation}
where $x_{1}^{\pm}(t),\ x_{2}^{\pm}(t)$ and $(\rho_*^{AB},v_*^{AB})$ can be determined by
(\ref{3.7}) and (\ref{3.8}).Let us use Fig.2(a) to illustrate this situation
in detail, where all the characteristics in the rarefaction wave
fans $R_1$ and $R_2$ are curved into
parabolic shapes.

If $(\rho_+,u_+)\in$ I\!I, the Riemann solutions to (\ref{1.1})-(\ref{1.2}) and
(\ref{1.6}) $S_1^{AB}+R_{2}^{AB}$ can be represented as
\begin{equation}\label{4.4}
(\rho,u)(x,t)=\left\{\begin{array}{ll}
(\rho_-,u_-+\beta t),\ \ \ \ \ \ \ \ \ \ \ x<x_{1}^{AB}(t),\\
(\rho_*,v_*+\beta t),\ \ \  \ \ \ \ \ \ \ \ \ x_{1}^{AB}(t)<x<x_{2}^-(t),\\
R_{2}^{AB},\ \ \  \ \ \ \ \ \ \ \ \ \ \ \ \ \  \ \ \ \ \ x_{2}^-(t)<x<x_{2}^+(t),\\
(\rho_+,u_++\beta t),\ \ \ \ \ \ \ \ \ \ \ x>x_{2}^+(t),
\end{array}\right.
\end{equation}
where $x_{1}^{AB}(t)$,$x_{2}^{\pm}(t)$ and $(\rho_*^{AB},v_*^{AB})$ are determined by (\ref{3.8}) and
(\ref{3.19}).
Let us use Fig.2(b) to illustrate this
situation in detail, where both the shock wave curve $S_1^{AB}$ and the characteristics in the rarefaction wave
fans $R_2^{AB}$ are curved into parabolic shapes.

If $(\rho_+,u_+)\in$ I\!I\!I, the Riemann solutions to (\ref{1.1})-(\ref{1.2})and
(\ref{1.6}) $R_1^{AB}+S_{2}^{AB}$, which is much similar to the case that $(\rho_+,u_+)\in$ I\!I, so we omit it, see Fig.2(c).

If $(\rho_+,u_+)\in$ I\!V, the Riemann solutions to (\ref{1.1})-(\ref{1.2}) and
(\ref{1.6}) $S_1^{AB}+S_{2}^{AB}$ can be represented as
\begin{equation}\label{4.2}
(\rho,u)(x,t)=\left\{\begin{array}{ll}
(\rho_-,u_-+\beta t),\ \ \ \ \ \ \ \ \ \ \ x<x_{1}^{AB}(t),\\
(\rho_*^{AB},v_*^{AB}+\beta t),\ \ \  \ \ \ \ \ x_{1}^{AB}(t)<x<x_{2}^{AB}(t),\\
(\rho_+,u_++\beta t),\ \ \ \ \ \ \ \ \ \ \ x>x_{2}^{AB}(t),
\end{array}\right.
\end{equation}
where $x_{1}^{AB}(t)$, $x_{2}^{AB}(t)$ and $(\rho_*^{AB},v_*^{AB})$ are determined by (\ref{3.18}),
(\ref{3.19}). Let us use Fig.2(d) to illustrate this
situation in detail, where both the shock wave curve $S_1^{AB}$ and the characteristics in the rarefaction wave
fans $S_2^{AB}$ are curved into parabolic shapes.


\unitlength 0.9mm 
\linethickness{0.4pt}
\ifx\plotpoint\undefined\newsavebox{\plotpoint}\fi 
\begin{picture}(174.5,102.5)(10,0)
\put(83,60.5){\vector(1,0){.07}}
\put(20.25,60.5){\line(1,0){62.75}}
\put(11.75,103){\vector(0,1){.07}}
\put(11.75,62.75){\line(0,1){40.25}}
\qbezier(47.75,60.5)(51.75,84.75)(77.25,94)
\qbezier(47.75,60.5)(60.25,81.63)(84,90.25)
\qbezier(47.75,60.5)(63.88,75.75)(88.5,85)
\qbezier(47.75,60.5)(26.88,64)(19.75,93)
\qbezier(47.75,60.5)(37.5,68.25)(31.25,95.5)
\qbezier(47.75,60.5)(31.75,66.75)(24.75,93.5)
\put(171.75,59.5){\vector(1,0){.07}}
\put(105.5,59.5){\line(1,0){66.25}}
\put(99.75,103){\vector(0,1){.07}}
\put(99.75,59.5){\line(0,1){42.5}}
\qbezier(132.5,59.5)(135.75,82.88)(152,96.5)
\qbezier(132.5,59.5)(142.38,80.5)(158.5,92)
\qbezier(132.5,59.5)(146.38,76.75)(161.75,84.75)
\qbezier(132.5,59.5)(119.13,66.63)(110.25,92.75)
\put(84.75,14.5){\vector(1,0){.07}}
\put(21,14.5){\line(1,0){63.75}}
\put(12.75,53.5){\vector(0,1){.07}}
\put(12.75,14.5){\line(0,1){37.75}}
\qbezier(48,14.5)(51.13,38.13)(73.75,42.5)
\qbezier(48,14.5)(39.75,19.75)(36.75,49)
\qbezier(48,14.5)(34.75,16.13)(29.5,47)
\qbezier(48,14.5)(24.38,16.13)(21.75,43.25)
\put(81.75,57.75){$x$}
\put(46.75,57){0}
\put(69,86){$R_2^{AB}$}
\put(28,86){$R_1^{AB}$}
\put(41,54){(a) $(\rho_+,u_+)\in$ I}
\put(8.75,100.75){$t$}
\put(36.25,86){$(\rho_\ast^{AB},u_\ast^{AB}+\beta t)$}
\put(12.75,64){$(\rho_-,u_-+\beta t)$}
\put(62,68.75){$(\rho_+,u_++\beta t)$}
\put(174,14.75){\vector(1,0){.07}}
\put(106,14.75){\line(1,0){68}}
\put(99.75,52){\vector(0,1){.07}}
\put(99.75,16){\line(0,1){36}}
\qbezier(137,14.75)(141.13,37.13)(161.75,44.5)
\qbezier(137,14.75)(121.38,24.13)(118.25,45.5)
\put(171.5,56){$x$}
\put(131.75,56){0}
\put(128.75,53){(b) $(\rho_+,u_+)\in$ I\!I}
\put(97,100){$t$}
\put(101,66){$(\rho_-,u_-+\beta t)$}
\put(108,86){$S_1^{AB}$}
\put(154,67.5){$(\rho_+,u_++\beta t)$}
\put(145,86){$R_2^{AB}$}
\put(118.75,81.5){$(\rho_\ast,u_\ast+\beta t)$}
\put(10.75,52.25){$t$}
\put(48,11.5){0}
\put(43,8){(c) $(\rho_+,u_+)\in$ I\!I\!I}
\put(83.75,12){$x$}
\put(13,17){$(\rho_-,u_-+\beta t)$}
\put(26,39){$R_2^{AB}$}
\put(39,36){$(\rho_\ast^{AB},u_\ast^{AB}+\beta t)$}
\put(52,25){$(\rho_+,u_++\beta t)$}
\put(63,41){$S_1^{AB}$}
\put(98.25,49.75){$t$}
\put(137.75,10.75){0}
\put(120,41){$S_1^{AB}$}
\put(160,41){$S_2^{AB}$}
\put(132.75,8){(d) $(\rho_+,u_+)\in$ I\!V}
\put(174.5,11.25){$x$}
\put(123,37.75){$(\rho_\ast^{AB},u_\ast^{AB}+\beta t)$}
\put(102,25){$(\rho_-,u_-+\beta t)$}
\put(145,27.75){$(\rho_+,u_++\beta t)$}
\put(50,3){Fig.2 The Riemann solutions of (\ref{1.1})-(\ref{1.2})
and (\ref{1.6}) when $\beta>0$.  }
\end{picture}


In a word, the Riemann problem for the extended Chaplygin gas equations with Coulomb-like friction is constructed completely.
 We see that the Riemann
solutions to (\ref{1.1})-(\ref{1.2}) and (\ref{1.6}) are not self-similar any
more, in which the state variable $u$ varies linearly along with the
time $t$ under the influence of the Coulomb-like friction. In other
words, the state variable $u-\beta t$ remains unchanged in the left,
intermediate and right states. It is shown that the Coulomb-like friction term
make shock waves and rarefaction waves bend into parabolic shapes for the Riemann
solutions.

\section{Formation of $\delta$-shocks and vacuum states as $A,B\rightarrow0$}

In this section, we will study the vanishing pressure limit process,
i.e.,$A,B\rightarrow0$. Since the two regions $\rm{I\!I}$ and $\rm{I\!I\!I}$ in the $(\rho,v)$ plane
have empty interior when $A,B\rightarrow0$, it suffices to analyze
the limit process for the two cases $(\rho_+,u_+)\in \rm{I}$ and $(\rho_+,u_+)\in\rm{I\!V}$.

\subsection{$\delta-$shocks and concentration}
Firstly, we analyze the formation of $\delta$-shocks in Riemann
solutions to the extended Chaplygin gas equations
$(\ref{1.1})$-$(\ref{1.2})$ with $(\ref{1.6})$ in the case
$(\rho_+,u_+)\in\rm{I\!V}$ and $u_->u_+$ as the pressure vanishes.

When $(\rho_+,u_+)\in \rm{I\!V} $, for fixed $A,B>0$,
let $(\rho_*^{AB},u_*^{AB})=(\rho_*^{AB},v_*^{AB}+\beta t)$ be the intermediate state in the sense
that $(\rho_-,u_-+\beta t)$ and $(\rho_*^{AB},v_*^{AB}+\beta t)$ are connected by
1-shock $S_1^{AB}$ with speed $\sigma_1^{AB}$, $(\rho_*^{AB},v_*^{AB}+\beta t)$
and $(\rho_+,u_++\beta t)$ are connected by 2-shock $S_2^{AB}$ with speed
$\sigma_2^{AB}$. Then it follows
\begin{equation}\label{5.1}S_1^{AB}:\ \
\left\{\begin{array}{ll}
\sigma_1^{AB}=\DF{\rho_{*}^{AB} v_{*}^{AB}-\rho_-u_-}{\rho_{*}^{AB}-\rho_-}+\beta t,\\
v_{*}^{AB}-u_-=-\sqrt{\frac{\rho_{*}^{AB}-\rho_{-}}{\rho_{*}^{AB}\rho_{-}}\Big(A((\rho_{*}^{AB})^{n}-\rho_{-}^{n})-
B(\frac{1}{(\rho_{*}^{AB})^{\alpha}}-\frac{1}{\rho_{-}^{\alpha}})\Big)},\
\ \rho_{*}^{AB}>\rho_-,
\end{array}\right.
\end{equation}

\begin{equation}\label{5.2}S_2^{AB}:\ \
\left\{\begin{array}{ll}
\sigma_2^{AB}=\DF{\rho_{+} u_{+}-\rho_{*}^{AB}v_{*}^{AB}}{\rho_{+}-\rho_{*}^{AB}}+\beta t,\\
u_{+}-v_{*}^{AB}=-\sqrt{\frac{\rho_{+}-\rho_{*}^{AB}}{\rho_{+}\rho_{*}^{AB}}\Big(A(\rho_{+}^{n}-(\rho_{*}^{AB})^{n})-
B(\frac{1}{\rho_{+}^{\alpha}}-\frac{1}{(\rho_{*}^{AB})^{\alpha}})\Big)},\
\ \rho_{+}<\rho_{*}^{AB}.
\end{array}\right.
\end{equation}

In the following, we give some lemmas to show the limit behavior of
the Riemann solutions of system $(\ref{1.1})$-$(\ref{1.2})$ with
$(\ref{1.6})$ as $A,B\rightarrow0$.

\begin{Lemma}\label{lem:5.1}
$\lim\limits_{A,B\rightarrow0}\rho_*^{AB}=+\infty.$
\end{Lemma}

\noindent\textbf{Proof.} Eliminating $v_{*}^{AB}$ in the second
equation of $(\ref{5.1})$ and $(\ref{5.2})$ gives
\begin{eqnarray}
u_{-}-u_{+}=&&\sqrt{\frac{\rho_{*}^{AB}-\rho_{-}}{\rho_{*}^{AB}\rho_{-}}
\Big(A((\rho_{*}^{AB})^{n}-\rho_{-}^{n})-B(\frac{1}{(\rho_{*}^{AB})^{\alpha}}-\frac{1}{\rho_{-}^{\alpha}})\Big)}\nonumber\\
&+&\sqrt{\frac{\rho_{+}-\rho_{*}^{AB}}{\rho_{+}\rho_{*}^{AB}}
\Big(A(\rho_{+}^{n}-(\rho_{*}^{AB})^{n})-B(\frac{1}{\rho_{+}^{\alpha}}-\frac{1}{(\rho_{*}^{AB})^{\alpha}})\Big)}.\label{5.3}
\end{eqnarray}
If
$\lim\limits_{A,B\rightarrow0}\rho_*^{AB}=M\in(\max\{\rho_{-},\rho_{+}\},+\infty)$,
then by taking the limit of $(\ref{5.3})$ as $A,B\rightarrow0$, we
obtain that $u_{-}-u_{+}=0$, which contradicts with $u_{-}>u_{+}$.
Therefore we must have
$\lim\limits_{A,B\rightarrow0}\rho_*^{AB}=+\infty.$

By Lemma \ref{lem:5.1}, from $(\ref{5.3})$ we immediately have the
following lemma.

\begin{Lemma}\label{lem:5.2}
$\lim\limits_{A,B\rightarrow0}A(\rho_*^{AB})^{n}=\DF{\rho_-\rho_+}{(\sqrt{\rho_-}+\sqrt{\rho_+})^2}(u_--u_+)^2$.
\end{Lemma}

\begin{Lemma}\label{lem:5.3}
Let $\frac{dx_1^{AB}(t)}{dt}=\sigma_1^{AB}(t)$,
$\frac{dx_2^{AB}(t)}{dt}=\sigma_2^{AB}(t)$, then
\begin{equation}\label{4.4}
\lim\limits_{A,B\rightarrow0}u_*^{AB}=\lim\limits_{A,B\rightarrow0}\sigma_1^{AB}=\lim\limits_{A,B\rightarrow0}\sigma_2^{AB}=\sigma_{0}+\beta t.
\end{equation}
\end{Lemma}

\noindent\textbf{Proof.} From the first equation of  $(\ref{5.1})$
and $(\ref{5.2})$ for $S_{1}^{AB}$ and $S_{2}^{AB}$, by Lemma \ref{lem:5.1},
we have
$$\lim\limits_{A,B\rightarrow0}\sigma_1^{AB}=
\lim\limits_{A,B\rightarrow0}\DF{\rho_{*}^{AB}
v_{*}^{AB}-\rho_-u_-}{\rho_{*}^{AB}-\rho_-}+\beta t=
\lim\limits_{A,B\rightarrow0}v_*^{AB}+\beta t=\lim\limits_{A,B\rightarrow0}u_*^{AB},$$

$$\lim\limits_{A,B\rightarrow0}\sigma_2^{AB}=
\lim\limits_{A,B\rightarrow0}\DF{\rho_+u_+-\rho_{*}^{AB}
v_{*}^{AB}}{\rho_+-\rho_{*}^{AB}}+\beta t=\lim\limits_{A,B\rightarrow0}v_*^{AB}+\beta t=\lim\limits_{A,B\rightarrow0}u_*^{AB},$$
which immediately lead to
$\lim\limits_{A,B\rightarrow0}u_*^{AB}=\lim\limits_{A,B\rightarrow0}\sigma_1^{AB}=\lim\limits_{A,B\rightarrow0}\sigma_2^{AB}$.

From the second equation of $(\ref{5.1})$, by Lemma
\ref{lem:5.1}-\ref{lem:5.2}, we get

\begin{eqnarray}
\lim\limits_{A,B\rightarrow0}v_*^{AB}&=&u_{-}-\lim\limits_{A,B\rightarrow0}\sqrt{\frac{\rho_{*}^{AB}-\rho_{-}}{\rho_{*}^{AB}\rho_{-}}
\Big(A((\rho_{*}^{AB})^{n}-\rho_{-}^{n})-B(\frac{1}{(\rho_{*}^{AB})^{\alpha}}-\frac{1}{\rho_{-}^{\alpha}})\Big)}\nonumber\\
&=&u_{-}-\sqrt{\frac{1}{\rho_{-}}\DF{\rho_-\rho_+}{(\sqrt{\rho_-}+\sqrt{\rho_+})^2}(u_--u_+)^2}\nonumber\\
&=&u_{-}-\DF{\sqrt{\rho_+}}{\sqrt{\rho_-}+\sqrt{\rho_+}}(u_--u_+)\nonumber\\
&=&\frac{\sqrt{\rho_-}u_-+\sqrt{\rho_+}u_+}{\sqrt{\rho_-}+\sqrt{\rho_+}}=\sigma_{0}\nonumber.
\end{eqnarray}
So $\lim\limits_{A,B\rightarrow0}u_*^{AB}=\lim\limits_{A,B\rightarrow0}v_*^{AB}+\beta t=\sigma_{0}+\beta t.$
The proof is completed.

\begin{Lemma}\label{lem:5.4}
\begin{equation}\label{5.5}
\lim\limits_{A,B\rightarrow0}\int_{x_1^{AB}}^{x_2^{AB}}\rho_*^{AB}dx=\sqrt{\rho_{+}\rho_{-}}(u_--u_+)t,
\end{equation}
\begin{equation}\label{5.6}
\lim\limits_{A,B\rightarrow0}\int_{x_1^{AB}}^{x_2^{AB}}\rho_*^{AB}u_*^{AB}dx=(\sigma_{0}+\beta t)\sqrt{\rho_{+}\rho_{-}}(u_--u_+)t.
\end{equation}

\end{Lemma}

\noindent\textbf{Proof.} The first equations of the Rankine-Hugoniot
condition $(\ref{3.11})$ for $S_{1}^{AB}$ and $S_{2}^{AB}$ read
\begin{equation}\nonumber
\left\{\begin{array}{ll}
\sigma_1^{AB}(\rho_*^{AB}-\rho_-)=\rho_{*}^{AB}
(v_{*}^{AB}+\beta t)-\rho_-(u_-+\beta t),\\
\sigma_2^{AB}(\rho_+-\rho_{*}^{AB})=\rho_+(u_++\beta t)-\rho_{*}^{AB}
(v_{*}^{AB}+\beta t),
\end{array}\right.
\end{equation}
from which we have
\begin{eqnarray}
\lim\limits_{A,B\rightarrow0}\rho_*^{AB}(\sigma_2^{AB}-\sigma_1^{AB})&=&\lim\limits_{A,B\rightarrow0}
(-\sigma_1^{AB}\rho_{-}+\sigma_2^{AB}\rho_{+}-\rho_+(u_++\beta t)+\rho_-(u_-+\beta t))\nonumber\\
&=&\sigma(\rho_{+}-\rho_{-})-(\rho_{+}u_+-\rho_{-}u_-)\nonumber\\
&=&\sqrt{\rho_{+}\rho_{-}}(u_--u_+).\nonumber
\end{eqnarray}
So
\begin{equation}\label{5.7}
\lim\limits_{A,B\rightarrow0}\int_{x_1^{AB}(t)}^{x_2^{AB}(t)}\rho_*^{AB}dx=
=\lim\limits_{A,B\rightarrow0}\int_{0}
^{t}\rho_*^{AB}(\sigma_2^{AB}-\sigma_1^{AB})dt=
\sqrt{\rho_{+}\rho_{-}}(u_--u_+)t.
\end{equation}
Similarly, from the second equations of the Rankine-Hugoniot
condition $(\ref{3.11})$ for $S_{1}^{AB}$ and $S_{2}^{AB}$, we have
\begin{equation}\nonumber
\left\{\begin{array}{ll}
\sigma_1^{AB}(\rho_*^{AB}v_{*}^{AB}-\rho_-u_-)=\rho_{*}^{AB}
v_{*}^{AB}(v_{*}^{AB}+\beta t)-\rho_-u_-(u_-+\beta t)+A((\rho_{*}^{AB})^{n}-\rho_-^{n})-B(\DF{1}{(\rho_{*}^{AB})^{\alpha}}-\DF{1}{\rho_-^{\alpha}}),\\
\sigma_2^{AB}(\rho_+u_+-\rho_{*}^{AB}v_{*}^{AB}v_{*}^{AB}=\rho_+u_+(u_++\beta t)-\rho_{*}^{AB}
v_{*}^{AB}(v_{*}^{AB}+\beta t)+A(\rho_+^{n}-(\rho_{*}^{AB})^{n})-B(\DF{1}{\rho_+^{\alpha}}-\DF{1}{(\rho_{*}^{AB})^{\alpha}}),
\end{array}\right.
\end{equation}

then we obtain
\begin{eqnarray}
&\lim\limits_{A,B\rightarrow0}&\rho_*^{AB}v_{*}^{AB}(\sigma_2^{AB}-\sigma_1^{AB})\nonumber\\
&=&\lim\limits_{A,B\rightarrow0}(-\sigma_1^{AB}\rho_{-}u_{-}+\sigma_2^{AB}\rho_{+}u_{+}-
\rho_+u_+(u_++\beta t)+\rho_-u_-(u_-+\beta t)-A(\rho_{+}^{n}-\rho_{-}^{n})+B(\DF{1}{\rho_+^{\alpha}}-\DF{1}{\rho_{-}^{\alpha}}))\nonumber\\
&=&\sigma_{0}(\rho_{+}u_+-\rho_{-}u_-)-(\rho_{+}u_+^2-\rho_{-}u_-^2)=\sigma_{0}\sqrt{\rho_{+}\rho_{-}}(u_--u_+).\nonumber
\end{eqnarray}
So
\begin{eqnarray}
\lim\limits_{A,B\rightarrow0}\int_{x_1^{AB}(t)}^{x_2^{AB}(t)}\rho_*^{AB}u_{*}^{AB}dx&=&\lim\limits_{A,B\rightarrow0}\int_{x_1^{AB}(t)}^{x_2^{AB}(t)}\rho_*^{AB}v_{*}^{AB}dx+
\beta t\lim\limits_{A,B\rightarrow0}\int_{x_1^{AB}(t)}^{x_2^{AB}(t)}\rho_*^{AB}dx\nonumber\\
&=&\lim\limits_{A,B\rightarrow0}\int_{0}
^{t}\rho_*^{AB}v_{*}^{AB}(\sigma_2^{AB}-\sigma_1^{AB})dt+\beta t\sqrt{\rho_{+}\rho_{-}}(u_--u_+)t\nonumber\\
&=&(\sigma_{0}+\beta t)\sqrt{\rho_{+}\rho_{-}}(u_--u_+)t.\label{5.8}
\end{eqnarray}
 The proof is finished.
\bigbreak

The above lemmas 5.1-5.4 show that, as $A,B\rightarrow0$, the curves of the
shock wave $S_{1}^{AB}$ and
$S_{2}^{AB}$ will coincide and the delta shock waves will form.  Next
we will arrange the values which gives the exact position,
propagation speed and strength of the delta shock wave according to
Lemma \ref{lem:5.3} and \ref{lem:5.4}.

 From (\ref{5.5}) and (\ref{5.6}), we let
\begin{equation}\label{5.9}
w(t)=\sqrt{\rho_{+}\rho_{-}}(u_--u_+)t,
\end{equation}
\begin{equation}\label{5.10}
w(t)u_\delta(t)=(\sigma_{0}+\beta t)\sqrt{\rho_{+}\rho_{-}}(u_--u_+)t,
\end{equation}
then
\begin{equation}\label{5.11}
u_\delta(t)=\sigma_{0}+\beta t,
\end{equation}
which is equal to $\sigma(t)$. Furthermore, by letting
$\frac{dx(t)}{dt}=\sigma(t)$, we have
\begin{equation}\label{5.12}
x(t)=\sigma_{0} t+\frac{1}{2}\beta t^2.
\end{equation}

From (\ref{5.9})-(\ref{5.12}), we can see that the quantities defined
above are exactly consistent with those given by
(\ref{2.7}).
Thus, it uniquely determines that the limits of the Riemann
solutions to the system (\ref{1.1})-(\ref{1.2}) and (\ref{1.6}) when
$A,B\rightarrow 0$ in the case $(\rho_+,u_+) \in$ I\!V and $u_{-}>u_{+}$ is just the
delta shock solution of (\ref{1.3}) and (\ref{1.6}). So we get the
following results which characterizes the vanishing pressure limit in
the case $(\rho_+,u_+)\in\rm{I\!V}$ and $u_{-}>u_{+}$.

\begin{thm}\label{5.1}
If $u_{-}>u_{+}$, for each fixed $A,B$, $(\rho_+,u_+)\in$
$\rm{I\!V}$, assuming that $(\rho^{AB},u^{AB})$ is a two-shock wave
solution of (\ref{1.1})-(\ref{1.2}) and (\ref{1.6}) which is constructed in
Section 4, it is obtained that when $A,B\rightarrow 0$, $(\rho^{AB},u^{AB})$
converges to a delta shock wave solution to the transportation
equations (\ref{1.3}) with the same source term and the same initial data.
\end{thm}

\subsection{Formation of vacuum states} In this subsection, we show the
formation of vacuum states in the Riemann solutions to
$(\ref{1.1})$-$(\ref{1.2})$ and $(\ref{1.6})$ in the case
$(\rho_+,u_+)\in I$ with $u_-<u_+$ and
$\rho_\pm>0$ as the pressure vanishes.

At this monent, for fixed $A,B>0$, let $(\rho_*^{AB},u_*^{AB})=(\rho_*^{AB},v_*^{AB}+\beta t)$ be
the intermediate state in the sense that $(\rho_-,u_-+\beta t)$ and
$(\rho_*^{AB},v_*^{AB}+\beta t)$ are connected by 1-rarefaction wave $R_1^{AB}$
with speed $\lambda_1^{AB}$, $(\rho_*^{AB},v_*^{AB}+\beta t)$and
$(\rho_+,u_++\beta t)$ are connected by 2-rarefaction wave $R_2^{AB}$ with speed
$\lambda_2^{AB}$. Then it follows
\begin{equation}\label{5.13}R_1^{AB}:\ \
\left\{\begin{array}{ll} \xi=\lambda_1^{AB}=v+\beta t-\sqrt{An\rho^{n-1}+\frac{\alpha B}{\rho^{\alpha+1}}},\\
v-u_-=-\int_{\rho_{-}}^{\rho}\frac{\sqrt{An\rho^{n-1}+\frac{\alpha
B}{\rho^{\alpha+1}}}}{\rho}d\rho,\ \
\rho_{*}^{AB}\leq\rho\leq\rho_{-}.
\end{array}\right.
\end{equation}

\begin{equation}\label{5.14}R_2^{AB}:\ \
\left\{\begin{array}{ll} \xi=\lambda_2^{AB}=v+\beta t+\sqrt{An\rho^{n-1}+\frac{\alpha B}{\rho^{\alpha+1}}},\\
u_+-v=\int_{\rho}^{\rho_{+}}\frac{\sqrt{An\rho^{n-1}+\frac{\alpha
B}{\rho^{\alpha+1}}}}{\rho}d\rho,\ \
\rho_{*}^{AB}\leq\rho\leq\rho_{+}.
\end{array}\right.
\end{equation}

Now, from the second equations of $(\ref{5.13})$ and $(\ref{5.14})$,
using the following integral identity
\begin{eqnarray}
&&\int^{\rho_{-}}_{\rho}\frac{\sqrt{An\rho_{-}^{n-1}+\frac{\alpha
B}{\rho^{\alpha+1}}}}{\rho}d\rho\nonumber\\
&=&\frac{2}{\alpha+1}\Big(-\sqrt{An\rho_{-}^{n-1}+\frac{\alpha
B}{\rho^{\alpha+1}}} +\sqrt{An\rho_{-}^{n-1}}
\ln(\sqrt{An\rho_{-}^{n-1}\rho^{\alpha+1}+\alpha
B}+\sqrt{An\rho_{-}^{n-1}\rho^{\alpha+1}})\Big)\Big
|_{\rho}^{\rho_{-}}\nonumber,
\end{eqnarray}

 it follows that the
intermediate state $(\rho_*^{AB},u_*^{AB})$ satisfies

\begin{eqnarray}
&&u_+-u_-\nonumber\\
&=&\int^{\rho_{-}}_{\rho_*^{AB}}\frac{\sqrt{An\rho^{n-1}+\frac{\alpha
B}{\rho^{\alpha+1}}}}{\rho}d\rho+\int_{\rho_*^{AB}}^{\rho_{+}}\frac{\sqrt{An\rho^{n-1}+\frac{\alpha
B}{\rho^{\alpha+1}}}}{\rho}d\rho\nonumber\\
&\leq&\int^{\rho_{-}}_{\rho_*^{AB}}\frac{\sqrt{An\rho_{-}^{n-1}+\frac{\alpha
B}{\rho^{\alpha+1}}}}{\rho}d\rho+\int_{\rho_*^{AB}}^{\rho_{+}}\frac{\sqrt{An\rho_{+}^{n-1}+\frac{\alpha
B}{\rho^{\alpha+1}}}}{\rho}d\rho\nonumber\\
&=&\frac{2}{\alpha+1}\Big(-\sqrt{An\rho_{-}^{n-1}+\frac{\alpha
B}{\rho_{-}^{\alpha+1}}} +\sqrt{An\rho_{-}^{n-1}}
\ln(\sqrt{An\rho_{-}^{n-1}\rho_{-}^{\alpha+1}+\alpha
B}+\sqrt{An\rho_{-}^{n-1}\rho_{-}^{\alpha+1}})\nonumber\\
&&+\sqrt{An\rho_{-}^{n-1}+\frac{\alpha B}{(\rho_*^{AB})^{\alpha+1}}}
-\sqrt{An\rho_{-}^{n-1}}
\ln(\sqrt{An\rho_{-}^{n-1}(\rho_*^{AB})^{\alpha+1}+\alpha
B}+\sqrt{An\rho_{-}^{n-1}(\rho_*^{AB})^{\alpha+1}})\nonumber\\
&&-\sqrt{An\rho_{+}^{n-1}+\frac{\alpha B}{\rho_{+}^{\alpha+1}}}
+\sqrt{An\rho_{+}^{n-1}}
\ln(\sqrt{An\rho_{+}^{n-1}\rho_{+}^{\alpha+1}+\alpha
B}+\sqrt{An\rho_{+}^{n-1}\rho_{+}^{\alpha+1}})\nonumber\\
&&+\sqrt{An\rho_{+}^{n-1}+\frac{\alpha B}{(\rho_*^{AB})^{\alpha+1}}}
-\sqrt{An\rho_{+}^{n-1}}
\ln(\sqrt{An\rho_{+}^{n-1}(\rho_*^{AB})^{\alpha+1}+\alpha
B}+\sqrt{An\rho_{+}^{n-1}(\rho_*^{AB})^{\alpha+1}})\Big
),\nonumber\\\label{5.15}
\end{eqnarray}
which implies that $$\lim\limits_{A,B\rightarrow0}\rho_*^{AB}=0$$.

Indeed, if
$\lim\limits_{A,B\rightarrow0}\rho_*^{AB}=K\in(0,\min\{\rho_{-},\rho_{+}\})$,
then $(\ref{5.15})$ leads to $u_+-u_-=0$, which contradicts with
$u_-<u_+$. Thus $\lim\limits_{A,B\rightarrow0}\rho_*^{AB}=0$, which
just means vacuum occurs. Moreover, as $A,B\rightarrow0$, one can
directly derive from $(\ref{5.13})$ and $(\ref{5.14})$ that
$\lambda_{1}^{AB},\ \lambda_{2}^{AB}\rightarrow u_{\pm}+\beta t$ and two
rarefaction waves $R_{1}^{AB}$ and $R_{2}^{AB}$ tend to two contact
discontinuities $\frac{dx}{dt}=u_\pm+\beta t$, respectively. These reach
the desired conclusion.

\begin{thm}\label{5.2}
Let $u_-<u_+$ and $(\rho_+,u_+)\in\rm{I}$. For any fixed
$A,B>0$, assume that $(\rho^{AB},u^{AB})$ is the two-rarefaction
wave Riemann solution of $(\ref{1.1})$-$(\ref{1.2})$ and $(\ref{1.6})$ with Riemann
data $(\rho_\pm,u_\pm)$ constructed in section 4. Then as
$A,B\rightarrow0$, the limit of the Riemann solution
 $(\rho^{AB}, u^{AB})$ is two contact
discontinuities connecting the states $(\rho_\pm,u_\pm+\beta t)$
and the intermediate vacuum state as follows:
\begin{equation}\label{2.1}
(\rho,u)(x,t)=\left\{\begin{array}{ll}
(\rho_-,u_-+\beta t),\ \ \ \ -\infty<x<u_-t+\frac{1}{2}\beta t^{2},\\
vacuum,\ \ \ \ \ \ \ \ \ \ u_-t+\frac{1}{2}\beta t^{2}<x<u_+t+\frac{1}{2}\beta t^{2},\\
(\rho_+,u_++\beta t),\ \ \ \ u_+t+\frac{1}{2}\beta t^{2}<x<\infty.
\end{array}\right.
\end{equation}
which is exactly the Riemann solution to the transportation equations
$(\ref{1.3})$ with the same source term and the same initial data.
\end{thm}

\section{Formation of $\delta$-shocks as $A\rightarrow0$}

In this section, we study the formation of the delta shock waves
in the limit of Riemann solutions of
$(\ref{1.1})$-$(\ref{1.2})$ and $(\ref{1.6})$ in the case
$(\rho_+,u_+)\in \rm{I\!V}$ and $  u_++\sqrt{\alpha
B}\rho_-^{-\frac{\alpha+1}{2}}<u_--\sqrt{\alpha
B}\rho_-^{-\frac{\alpha+1}{2}}$ as the pressure
approaches the generalized Chaplygin gas pressure, i.e.,
$A\rightarrow0$.

When $(\rho_+,u_+) \in\rm{I\!V} $, for fixed $A,B>0$,
let $(\rho_*^{A},u_*^{A})=(\rho_*^{A},v_*^{A}+\beta t)$ be the intermediate state in the sense
that $(\rho_-,u_-+\beta t)$ and $(\rho_*^{A},v_*^{A}+\beta t)$ are connected by
1-shock $S_1^{A}$ with speed $\sigma_1^{A}$, $(\rho_*^{A},v_*^{A}+\beta t)$
and $(\rho_+,u_++\beta t)$ are connected by 2-shock $S_2^{A}$ with speed
$\sigma_2^{A}$. Then it follows

\begin{equation}\label{6.1}S_1^{A}:\ \
\left\{\begin{array}{ll}
\sigma_1^{A}=\DF{\rho_{*}^{A} v_{*}^{A}-\rho_-u_-}{\rho_{*}^{A}-\rho_-}+\beta t,\\
v_{*}^{A}-u_-=-\sqrt{\frac{\rho_{*}^{A}-\rho_{-}}{\rho_{*}^{A}\rho_{-}}\Big(A((\rho_{*}^{A})^{n}-\rho_{-}^{n})-
B(\frac{1}{(\rho_{*}^{A})^{\alpha}}-\frac{1}{\rho_{-}^{\alpha}})\Big)},\
\ \rho_{*}^{A}>\rho_-,
\end{array}\right.
\end{equation}
and
\begin{equation}\label{6.2}S_2^{A}:\ \
\left\{\begin{array}{ll}
\sigma_2^{A}=\DF{\rho_{+} u_{+}-\rho_{*}^{A}v_{*}^{A}}{\rho_{+}-\rho_{*}^{A}}+\beta t,\\
u_{+}-v_{*}^{A}=-\sqrt{\frac{\rho_{+}-\rho_{*}^{A}}{\rho_{+}\rho_{*}^{A}}\Big(A(\rho_{+}^{n}-(\rho_{*}^{A})^{n})-
B(\frac{1}{\rho_{+}^{\alpha}}-\frac{1}{(\rho_{*}^{A})^{\alpha}})\Big)},\
\ \rho_{+}<\rho_{*}^{A}.
\end{array}\right.
\end{equation}
Here $\sigma_1^{A}$ and $\sigma_2^{A}$ are the propagation speed of
$S_{1}^{A}$ and $S_{2}^{A}$, respectively. Similar to that in Section 5,
 in the following, we give some lemmas to show the limit behavior
of the Riemann solutions of system $(\ref{1.1})$-$(\ref{1.2})$ and
$(\ref{1.6})$ as $A\rightarrow0$.

\begin{Lemma}\label{lem:6.1}
$\lim\limits_{A\rightarrow0}\rho_*^{A}=+\infty.$
\end{Lemma}

\noindent\textbf{Proof.} Eliminating $v_{*}^{A}$ in the second
equation of $(\ref{6.1})$ and $(\ref{6.2})$ gives
\begin{eqnarray}
u_{-}-u_{+}=&&\sqrt{\frac{\rho_{*}^{A}-\rho_{-}}{\rho_{*}^{A}\rho_{-}}
\Big(A((\rho_{*}^{A})^{n}-\rho_{-}^{n})-B(\frac{1}{(\rho_{*}^{A})^{\alpha}}-\frac{1}{\rho_{-}^{\alpha}})\Big)}\nonumber\\
&+&\sqrt{\frac{\rho_{+}-\rho_{*}^{A}}{\rho_{+}\rho_{*}^{A}}
\Big(A(\rho_{+}^{n}-(\rho_{*}^{A})^{n})-B(\frac{1}{\rho_{+}^{\alpha}}-\frac{1}{(\rho_{*}^{A})^{\alpha}})\Big)}.\label{6.3}
\end{eqnarray}
If
$\lim\limits_{A\rightarrow0}\rho_*^{A}=K\in(\max\{\rho_{-},\rho_{+}\},+\infty)$,
then by taking the limit of $(\ref{6.3})$
\begin{eqnarray}
u_{-}-u_{+}&=&\sqrt{B}\Big(\sqrt{(\frac{1}{\rho_{-}}-\frac{1}{K})(\frac{1}{\rho_{-}^{\alpha}}-\frac{1}{K^{\alpha}})}+
\sqrt{(\frac{1}{\rho_{+}}-\frac{1}{K})(\frac{1}{\rho_{+}^{\alpha}}-\frac{1}{K^{\alpha}})}\Big)\nonumber\\
&<&\sqrt{B}\Big(\sqrt{\frac{1}{\rho_{-}}\frac{1}{\rho_{-}^{\alpha}}}+
\sqrt{\frac{1}{\rho_{+}}\frac{1}{\rho_{+}^{\alpha}}}\Big)\nonumber\\
&=&\sqrt{ B}\rho_+^{-\frac{\alpha+1}{2}} +\sqrt{
B}\rho_-^{-\frac{\alpha+1}{2}}\label{6.4}
\end{eqnarray}
which contradicts with the condition $u_++\sqrt{\alpha
B}\rho_-^{-\frac{\alpha+1}{2}}<u_--\sqrt{\alpha
B}\rho_-^{-\frac{\alpha+1}{2}}$. Therefore we must have
$\lim\limits_{A\rightarrow0}\rho_*^{A}=+\infty.$

By Lemma \ref{lem:6.1}, from $(\ref{6.3})$ we immediately have the
following lemma.

\begin{Lemma}\label{lem:6.2}
$\lim\limits_{A\rightarrow0}A(\rho_*^{A})^{n}<\rho_-(u_--u_+)^2$.
\end{Lemma}

\begin{Lemma}\label{lem:6.3}
Let $ \frac{dx_{1}^{A}}{dt}=\sigma_1^{A}, \ \frac{dx_{2}^{A}}{dt}=\sigma_2^{A}, \ \lim\limits_{A\rightarrow0}v_*^{A}=\widehat{\sigma_0^{B}}$, then
\begin{equation}\label{6.5}
\lim\limits_{A\rightarrow0}u_*^{A}=\lim\limits_{A\rightarrow0}\sigma_1^{A}
=\lim\limits_{A\rightarrow0}\sigma_2^{A}=\widehat{\sigma_{0}^{B}}+\beta t,
\end{equation}
\end{Lemma}
where $\widehat{\sigma_{0}^{B}}\in\Big(u_++\sqrt{
\alpha B}\rho_+^{-\frac{\alpha+1}{2}}, u_--\sqrt{\alpha
B}\rho_-^{-\frac{\alpha+1}{2}}\Big)$.

\noindent\textbf{Proof.} From the second equation of $(\ref{6.1})$
for $S_{1}^{A}$, by Lemma \ref{lem:6.1} and \ref{lem:6.2}, we have
\begin{eqnarray}
\lim\limits_{A\rightarrow0}v_*^{A}&=&u_{-}-\lim\limits_{A\rightarrow0}\sqrt{\frac{\rho_{*}^{A}-\rho_{-}}{\rho_{*}^{A}\rho_{-}}
\Big(A((\rho_{*}^{A})^{n}-\rho_{-}^{n})-B(\frac{1}{(\rho_{*}^{A})^{\alpha}}-\frac{1}{\rho_{-}^{\alpha}})\Big)}\nonumber\\
&=&u_{-}-\sqrt{\frac{1}{\rho_{-}}\Big(\lim\limits_{A\rightarrow0}A(\rho_*^{A})^{n}+\frac{B}{\rho_{-}^{\alpha}}\Big)}\nonumber\\
&<& u_--\sqrt{\alpha
B}\rho_-^{-\frac{\alpha+1}{2}}\nonumber\\\label{6.6}.
\end{eqnarray}

Similarly, from the second equation of $(\ref{6.2})$ for $S_{2}^{A}$, we
have
\begin{eqnarray}
\lim\limits_{A\rightarrow0}v_*^{A}&=&u_{+}+\lim\limits_{A\rightarrow0}\sqrt{\frac{\rho_{+}-\rho_{*}^{A}}{\rho_{+}\rho_{*}^{A}}\Big(A(\rho_{+}^{n}-(\rho_{*}^{A})^{n})-
B(\frac{1}{\rho_{+}^{\alpha}}-\frac{1}{(\rho_{*}^{A})^{\alpha}})\Big)}\nonumber\\
&=&u_{+}+\sqrt{\frac{1}{\rho_{+}}\Big(\lim\limits_{A\rightarrow0}A(\rho_*^{A})^{n}+\frac{B}{\rho_{+}^{\alpha}}\Big)}\nonumber\\
&>&u_++\sqrt{ \alpha B}\rho_+^{-\frac{\alpha+1}{2}}
\nonumber\\\label{6.7}.
\end{eqnarray}
Furthermore, similar to the analysis as Lemma \ref{lem:5.3}, we can
obtain
$\lim\limits_{A\rightarrow0}u_*^{A}=\lim\limits_{A\rightarrow0}\sigma_1^{A}=\lim\limits_{A\rightarrow0}\sigma_2^{A}=\widehat{\sigma_{0}^{B}}+\beta t$.
The proof is complete.

\begin{Lemma}\label{lem:6.4}
For $\widehat{\sigma_{0}^{B}}$ mentioned in Lemma \ref{lem:6.3},
\begin{equation}\label{6.8}
\widehat{\sigma_{0}^{B}}=\sigma_{0}^{B}=\DF{\rho_+ u_+-\rho_-u_-+
\big\{\rho_+\rho_-\big((u_+-u_-)^2-(\frac{1}{\rho_+}-\frac{1}{\rho_-})
(\frac{B}{\rho_+^\alpha}-\frac{B}{\rho_-^\alpha})\big)\big\}^{\frac{1}{2}}}{\rho_+-\rho_-},
\end{equation}
as $\rho_+\neq\rho_-$, and
\begin{equation}\label{6.9}
\widehat{\sigma_{0}^{B}}=\sigma_{0}^{B}=\frac{u_++u_-}{2}
\end{equation}
as $\rho_+=\rho_-$.
\end{Lemma}

\noindent\textbf{Proof.} Let
$\lim\limits_{A\rightarrow0}A(\rho_*^{A})^{n}=L$, by Lemma
\ref{lem:6.3}, from $(\ref{6.6})$ and $(\ref{6.7})$ we have
$$\lim\limits_{A\rightarrow0}v_*^{A}=
u_{-}-\sqrt{\frac{1}{\rho_{-}}\Big(L+\frac{B}{\rho_{-}^{\alpha}}\Big)}=
u_{+}+\sqrt{\frac{1}{\rho_{+}}\Big(L+\frac{B}{\rho_{+}^{\alpha}}\Big)}=\widehat{\sigma_{0}^{B}},$$
which leads to
\begin{equation}\label{6.10}
L+\frac{B}{\rho_{+}^{\alpha}}=\rho_{-}(u_--\widehat{\sigma_{0}^{B}})^{2},
\end{equation}
\begin{equation}\label{6.11}
L+\frac{B}{\rho_{-}^{\alpha}}=\rho_{+}(u_+-\widehat{\sigma_{0}^{B}})^{2}.
\end{equation}
Eliminating $L$ from $(\ref{6.10})$ and $(\ref{6.11})$, we have
\begin{equation}\label{6.12}
(\rho_+-\rho_-)(\widehat{\sigma^{B}})^{2}-2(\rho_+
u_+-\rho_-u_-)\widehat{\sigma^{B}}+\rho_+
u_+^{2}-\rho_-u_-^{2}-B(\frac{1}{\rho_+^\alpha}-\frac{1}{\rho_-^\alpha})=0.
\end{equation}

From $(\ref{6.12})$, noticing
$\widehat{\sigma_{0}^{B}}\in\Big(u_++\sqrt{ \alpha
B}\rho_+^{-\frac{\alpha+1}{2}}, u_--\sqrt{\alpha
B}\rho_-^{-\frac{\alpha+1}{2}}\Big)$, we immediately get
$(\ref{6.8})$ and $(\ref{6.9})$. The proof is finished.

\bigbreak
Similar to  Lemma \ref{lem:5.4}, we have the following lemma.
\begin{Lemma}\label{lem:6.5}
\begin{equation}\label{6.13}
\lim\limits_{A\rightarrow0}\int_{x_1^{A}}^{x_2^{A}}\rho_*^{A}dx=w_0^{B}t,
\end{equation}
\begin{equation}\label{6.14}
\lim\limits_{A\rightarrow0}\int_{x_1^{A}}^{x_2^{A}}\rho_*^{A}u_*^{A}dx=(\sigma_{0}^{B}+\beta t)w_0^{B}t,
\end{equation}
\end{Lemma}

\noindent\textbf{Proof.} Here we only prove the case for $\rho_+\neq\rho_-$.
Similar to the proof of Lemma \ref{lem:5.4}, taking account into $(\ref{2.15})$ and $(\ref{6.12})$, we have
\begin{eqnarray}
\lim\limits_{A\rightarrow0}\rho_*^{A}(\sigma_2^{A}-\sigma_1^{A})&=&\lim\limits_{A\rightarrow0}
(-\sigma_1^{A}\rho_{-}+\sigma_2^{A}\rho_{+}-\rho_+(u_++\beta t)+\rho_-(u_-+\beta t))\nonumber\\
&=&\sigma_{0}^{B}(\rho_{+}-\rho_{-})-(\rho_{+}u_+-\rho_{-}u_-)=w_{0}^B\nonumber,
\end{eqnarray}
and
\begin{eqnarray}
&\lim\limits_{A\rightarrow0}&\rho_*^{A}v_{*}^{A}(\sigma_2^{A}-\sigma_1^{A})\nonumber\\
&=&\lim\limits_{A\rightarrow0}(-\sigma_1^{A}\rho_{-}u_{-}+\sigma_2^{A}\rho_{+}u_{+}-
\rho_+u_+(u_++\beta t)+\rho_-u_-(u_-+\beta t)-A(\rho_{+}^{n}-\rho_{-}^{n})+B(\DF{1}{\rho_+^{\alpha}}-\DF{1}{\rho_{-}^{\alpha}}))\nonumber\\
&=&\sigma_{0}^{B}(\rho_{+}u_+-\rho_{-}u_-)-(\rho_{+}u_+^2-\rho_{-}u_-^2)+B(\DF{1}{\rho_+^{\alpha}}-\DF{1}{\rho_{-}^{\alpha}})\nonumber\\
&=&(\rho_+-\rho_-)(\sigma_{0}^{B})^{2}-(\rho_+
u_+-\rho_-u_-)\sigma_{0}^{B}\nonumber\\
&=&\sigma_{0}^{B}(\sigma_{0}^{B}(\rho_+-\rho_-)-(\rho_+
u_+-\rho_-u_-))\nonumber\\
&=&\sigma_{0}^{B}w_{0}^{B}\nonumber.
\end{eqnarray}
So
\begin{equation}\label{6.15}
\lim\limits_{A\rightarrow0}\int_{x_1^{A}(t)}^{x_2^{A}(t)}\rho_*^{A}dx=
=\lim\limits_{A\rightarrow0}\int_{0}
^{t}\rho_*^{A}(\sigma_2^{A}-\sigma_1^{A})dt=
w_{0}^{B}t.
\end{equation}

\begin{eqnarray}
\lim\limits_{A\rightarrow0}\int_{x_1^{A}(t)}^{x_2^{A}(t)}\rho_*^{A}u_{*}^{A}dx&=&\lim\limits_{A\rightarrow0}\int_{x_1^{A}(t)}^{x_2^{A}(t)}\rho_*^{A}v_{*}^{A}dx+
\beta t\lim\limits_{A\rightarrow0}\int_{x_1^{A}(t)}^{x_2^{A}(t)}\rho_*^{A}dx\nonumber\\
&=&\lim\limits_{A\rightarrow0}\int_{0}
^{t}\rho_*^{A}v_{*}^{A}(\sigma_2^{A}-\sigma_1^{A})dt+\beta tw_{0}^{B}t\nonumber\\
&=&\sigma_{0}^{B}w_{0}^{B}t+\beta tw_{0}^{B}t\nonumber\\
&=&(\sigma_{0}^{B}+\beta t)w_{0}^{B}t.\label{6.16}
\end{eqnarray}
For the case $\rho_+=\rho_-$, the conclusion is obviously true, so we omit it. The proof is finished.
\bigbreak

The above lemmas 6.1-6.5 show that, as $A\rightarrow0$, the curves of the
shock wave $S_{1}^{A}$ and
$S_{2}^{A}$ will coincide and the delta shock waves will form.  Next,
we will arrange the values which give the exact position,
propagation speed and strength of the delta shock wave according to
Lemma \ref{lem:6.3} and \ref{lem:6.5}.

 From (\ref{6.13}) and (\ref{6.14}), when $\rho_+\neq\rho_-$, we let
\begin{equation}\label{6.17}
w^B(t)=w_{0}^{B}t,
\end{equation}
\begin{equation}\label{6.18}
w^B(t)u_\delta^B(t)=(\sigma_{0}^{B}+\beta t)w_{0}^{B}t,
\end{equation}
then
\begin{equation}\label{6.19}
u_\delta^B(t)=\sigma_{0}^B+\beta t,
\end{equation}
which is equal to $\sigma^B(t)$. Furthermore, by letting
$\frac{dx^B(t)}{dt}=\sigma^B(t)$, we have
\begin{equation}\label{6.20}
x^B(t)=\sigma_{0}^B t+\frac{1}{2}\beta t^2.
\end{equation}

From (\ref{6.17})-(\ref{6.20}), we can see that the quantities defined
above are exactly consistent with those given by (\ref{2.11})-
(\ref{2.16}). When $\rho_+=\rho_-$, similar result can be obtained.
Thus, it uniquely determines that the limits of Riemann
solutions to the system (\ref{1.1})-(\ref{1.2}) and (\ref{1.6}) when
$A\rightarrow 0$ in the case $(\rho_+,u_+)\in$ I\!V and $u_{-}> u_{+}$ is just the
delta shock solution of (\ref{1.5}) and (\ref{1.6}). So we get the
following results which characterizes the vanishing pressure limit in
the case $(\rho_+,u_+)\in\rm{I\!V}$ and $u_{-}> u_{+}$.

\begin{thm}\label{6.1}
If $u_{-}> u_{+}$, for each fixed $A,B$, $(\rho_+,u_+)\in$
$\rm{I\!V}$, assuming that $(\rho^{A},u^{A})$ is a two-shock wave
solution of (\ref{1.1})-(\ref{1.2}) and (\ref{1.6}) which is constructed in
Section 4, it is obtained that when $A\rightarrow 0$, $(\rho^{A},u^{A})$
converges to a delta shock wave solution to the generalized Chaplygin gas
equations (\ref{1.5}) with the same source term and the same initial data.
\end{thm}

\section{Conclusions and discussions}
In this paper, we have not only constructed the Riemann solutions for the nonhomogeneous
extended Chaplygin gas equations with coulomb-like friction, but also
considered two kinds of the flux approximation
limit of Riemann solutions to extended Chaplygin gas equations with coulomb-like friction and
studied the concentration and the formation of delta shocks during
the limit process. Moreover, we have proved that the vanishing
pressure limit of Riemann solutions to extended Chaplygin gas
equations with coulomb-like friction is just the corresponding ones to transportation equations with coulomb-like friction, and
when the extended Chaplygin pressure approaches the generalized
Chaplygin pressure, the limit of the two-shock Riemann solutions to the extended
Chaplygin gas equations with coulomb-like friction is just the delta shock wave solution to the
generalized Chaplygin gas equations with coulomb-like friction. In a word, in this paper, we have generalized all the
results about the vanishing pressure limit now available for homogeneous equations to the nonhomogeneous case, which is just the novelty of this paper
lies in.

In fact, one can further prove
that when the extended Chaplygin pressure approaches the pressure
for the perfect fluid, i.e., $B\rightarrow0$ for fixed $A$, the
limit of Riemann solutions to the nonhomogeneous extended Chaplygin gas
equations is just the corresponding ones to the nonhomogeneous Euler equations for
perfect fluids. Furthermore, the method of this paper can also be applied to other conservation laws with source term, such as transportation eqauations or shallow water equations with source term, which we will study in the future.

%

 \bigbreak

\noindent{\bf Acknowledgements}

The author is grateful to the anonymous referees for his/her valuable comments and corrections, which
helped to improve the manuscript. This work is partially supported by National Natural
Science Foundation of China (71601085).

\bigbreak


\begin{thebibliography}{99}


\bibitem{Benaoum}H.B. Benaoum. Accelerated universe from modified Chaplygin gas and
tachyonic fluid. {\it arXiv: hep-th/0205140}.

\bibitem{Bilic-Tupper-Viollier} N. Bilic, G.B. Tupper and R. Viollier. Dark matter, dark energy and the Chaplygin gas. {\it arXiv: astro-ph/0207423}.






\bibitem{Bouchut} F. Bouchut. On zero-pressure gas dynamics, in: Advances in Kinetic Theory and Computing.
{\it Ser. Adv. Math. Appl. Sci.} {\bf 22}, World Sci. Publishing:
River Edge, NJ, (1994), 171-190.

\bibitem{Brenier} Y. Brenier. Solutions with concentration to the Riemann problem for one-dimensional Chaplygin gas equations {\it J.
Math. Fluid Mech.} {\bf 7} (2005), S326-S331.

\bibitem{Brenier-Grenier} Y. Brenier and E. Grenier. Sticky particles and scalar conservation laws.
{\it SIAM J. Numer. Anal.} {\bf 35} (1998), 2317-2328.

\bibitem{Chang-Chen-Yang1} T. Chang, G.Q. Chen and S. Yang. On the Riemann problem for two-dimensional Euler equations. I. Interaction of shocks
and rarefaction waves. {\it Discrete Contin. Dyn. Syst.} {\bf 1}
(1995), 555-584.

\bibitem{Chang-Chen-Yang2} T. Chang, G.Q. Chen and S. Yang. On the Riemann problem for two-dimensional Euler equations. II. Interaction of contact
discontinuties. {\it Discrete Contin. Dyn. Syst.} {\bf 6} (2000),
419-430.

\bibitem{Chaplygin} S. Chaplygin. On gas jets. {\it Sci. Mem. Moscow Univ. Math. Phys.} {\bf 21} (1904), 1-121.

\bibitem{Chen-Liu1} G.Q. Chen and H. Liu. Formation of $\delta-$shocks and vacuum states in the vanishing pressure limit of solutions to
the Euler equations for isentropic fluids. {\it SIAM J. Math. Anal.}
{\bf 34} (2003), 925-938.

\bibitem{Chen-Liu2} G.Q. Chen and H. Liu. Concentration and cavitation in the vanishing pressure limit of solutions to
the Euler equations for nonisentropic fluids. {\it Physica D} {\bf
189} (2004), 141-165.

\bibitem{Chen-Sheng} J. Chen and W.C.Sheng. The Riemann problem and the limit solutions
as magnetic field vanishes to magnetogasdynamics for generalized Chaplygin gas. {\it Commun. Pure Appl. Anal.} {\bf
1} (2018), 127-142.

\bibitem{Danilvo-Shelkovich1}V.G. Danilvo and V.M. Shelkovich. Dynamics of propagation and interaction of
$\delta$-shock waves in conservation law system. {\it J.
Differential Equations} {\bf 221} (2005), 333-381.

\bibitem{Danilvo-Shelkovich2} V.G. Danilov and V.M. Shelkovich,
   Delta-shock waves type solution of hyperbolic systems of conservation
   laws, {\it Q.Appl.Math. }, 63(2005), 401-427.



\bibitem{E-Rykov-Sinai} W. E, Yu.G. Rykov and Ya.G. Sinai.
Generalized varinational principles, global weak solutions and
behavior with random initial data for systems of conservation laws
arising in adhesion particle dynamics. {\it Comm. Math. Phys. }{\bf
177} (1996), 349-380.

\bibitem{Faccanoni-Mangeney} G.Faccanoni, A.Mangeney. Exact solution for granular flows, {\it Int. J. Numer. Anal. Mech.
Geomech.}, {\bf 37} (2012), 1408-1433.

\bibitem{Gorini-Kamenshchik-Moschella-Pasquier} V. Gorini, A. Kamenshchik, U. Moschella and V.Pasquier.
The Chaplygin gas as a model for dark energy. {\it arXiv:
gr-qc/0403062}.

\bibitem{Guo-Li-Yin1} L.Guo, T.Li and G.Yin. The vanishling pressure limits of Riemann solutions
to the Chaplygin gas equations with a source term. {\it Commun. Pure Appl. Anal.} {\bf
16} (2017), 295-309.

\bibitem{Guo-Li-Yin2} L.Guo, T.Li and G.Yin. The limit behavior of the Riemann solutions
to the genralized Chaplygin gas equations with a source term. {\it J. Math. Anal. Appl.} {\bf
455} (2017), 127-140.


\bibitem{Guo-Sheng-Zhang} L.H. Guo, W.C. Sheng and T. Zhang. The two-dimensional Riemann problelm for
isentropic Chaplygin gas dynamic system. {\it Comm. Pure Appl.
Anal.} {\bf 9} (2) (2010), 431-458.




\bibitem{Huang-Wang} F. Huang and Z. Wang. Well posedness for
pressureless flow. {\it Comm. Math. Phys.} {\bf 222} (2001),
117-146.

\bibitem{von Karman}T.von Karman. Compressibility effects in aerodynamics {\it J. Aeron.
Sci.}, {\bf 8}(1941), pp.337-365.

\bibitem{Korchinski} D.J. Korchinski. Solutions of a Riemann problem for a system of conservation
laws possessing no classical weak solution. {\it Thsis, Adelphi
University,} (1977).

\bibitem{LeFloch-Liu}P.G. LeFloch and T.P. Liu. Existence theory to nonlinear hyperlolic
systems under nonconservative form. {\it Forum Math.} {\bf 5}
(1993), 261-280.

%

\bibitem{Li-Cao} Y. Li and Y. Cao. Large partial difference method with second accuracy in gas dynamics.
{\it Sci. Sinica A} {\bf 28} (1985), 1024-1035.


\bibitem{Naji}J.Naji. Extended Chaplygin gas equation of state with bulk and shear viscosities. {\it Astrophys. Space Sci.} {\bf 350} (2014),
333-338.

\bibitem{Qu-Wang}A.F. Qu and Z. Wang. Stability of the Riemann solutions for a
Chaplygin gas.  {\it J. Math. Anal. Appl.}, {\bf 409} (2014),
347-361.

\bibitem{Setare}M.R. Setare. Interacting holographic generalized Chaplygin gas
model. {\it Phys. Lett. B} {\bf 654} (2007), 1-6.




\bibitem{Shelkovich} V.M. Shelkovich. $\delta-$ and $\delta^{'}-$
shock wave types of singular solutions of systems of conservation
laws and transport and concentration processes. {\it Russian Math.
Surveys}, {\bf 63} (2008),473-546.

\bibitem{Shen1} C.Shen. The Riemann problem for the pressureless Euler system with the Coulomb-like friction term.
{\it IAM J. Appl. Math.}, 2015:1-24,doi:10.1039/imamat/hxv028.



\bibitem{Shen2} C.Shen. The Riemann problem for the Chaplygin gas equations with a source term.
{\it Z. Angew. Math. Mech.}, 2015:1-15,doi:10.1002/zamm.201500015.


\bibitem{Shen-Sun}C. Shen and M.N. Sun. Formation of delta shocks and vacuum states in
the vanishing pressure limit of Riemann solutions to the perturbed
Aw-Rascle model. {\it J. Differential Equations}, {\bf 249} (2010),
3024-3051.


\bibitem{Sheng-Zhang} W. Sheng and T. Zhang. The Riemann problem for transportation equations in gas dynamics.
{\it Mem. Amer. Math. Soc.} {\bf 137} (654) (1999), AMS: Providence.

\bibitem{Smoller} J. Smoller.
{\it Shock Waves and Reaction-Diffusion Equation.} Springer-Verlag,
New York, 1994.

\bibitem{Sun} M.Sun. The exact Riemann solutions to the generalized Chaplygin gas equations with friction.
{\it Commun. Nonlinear Sci. Numer. Simulat.}, 2016, {\bf 36}:
342-353.

\bibitem{Tan-Zhang-Zheng} D. Tan, T. Zhang and Y. Zheng. Delta-shock wave as limits of vanishing
viscosity for hyperbolic system of conservation laws. {\it J.
Differential Equations} {\bf 112} (1) (1994), 1-32.

\bibitem{Tsien} H.S. Tsien. Two dimensional subsonic flow of compressible fluids.
{\it J. Aeron. Sci.} {\bf 6} (1939), 399-407.

\bibitem{Wang}G. Wang. The Riemann problem for one dimensional
generalized Chaplygin gas dynamics. {\it J. Math. Anal. Appl.} {\bf
403} (2)(2013), 434-450.


\bibitem{Wang-Ding} Z. Wang and X. Ding. Uniqueness of generalized
solution for the Cauchy problem of transportation equations. {\it
Acta Math. Scientia} {\bf 17} (3) (1997), 341-352.

\bibitem{Wang-Huang-Ding} Z. Wang, F. Huang and X. Ding. On the
Cauchy problem of transportation equations. {\it Acta Math. Appl.
Sinica} {\bf 13} (2) (1997), 113-122.




\bibitem{Wang-Zhang2}Z. Wang and Q.L. Zhang. The Riemann problem with delta initial data
for the one-dimensional Chaplygin gas equations, {\it Acta
Mathematica Scientia}, {\bf 32B(3)} (2012):pp.825-841.

\bibitem{Yang-Wang}H. Yang and J. Wang. Delta -shocks and vacuum states in the vanishing
pressure limit of solutions to the isentropic Euler equations for
modified Chaplygin gas. {\it J. Math. Anal. Appl.} {\bf 413} (2014),
800-820.

\bibitem{Yang-Wang1}H. Yang and J. Wang. Concentration in the vanishing
pressure limit of solutions to the modified Chaplygin gas equations.
{\it J. Math. Phys.} {\bf 57} (111504)(2016); doi:10.1063/1.4967299.

\bibitem{Yin-Sheng} G. Yin and W. Sheng. $\delta-$shocks and vacuum states in the vanishing pressure limits of solutions to
the relativistic Euler equations for polytropic gases. {\it J. Math.
Anal. Appl.} {\bf 355} (2009), 596-605.




\end{thebibliography}
\end{document}